\documentclass{article}

\usepackage{arxiv}

\usepackage{amssymb}
\usepackage{amsthm}

\usepackage{float}
\usepackage{amsmath}
\usepackage{algorithm,algorithmic}
\usepackage{xcolor}
\usepackage{epsfig}
\usepackage{latexsym,bm,multirow}
\usepackage{doi}

\graphicspath{{./fig/}}

\DeclareMathOperator{\sech}{sech}
\newcommand{\bbR}{\mathbb{R}}
\newcommand{\rmd}{\mathrm{d}}

\usepackage{hyperref}
\title{Structure-preserving reduced-order modelling of Korteweg-de Vries equation}

\author{Murat Uzunca\\
   Department of Mathematics, Sinop University\\
     Sinop-Turkey\\ \
     \texttt{muzunca@sinop.edu.tr}\\
	\And
 S\"uleyman Y{\i}ld{\i}z\\
Institute of Applied Mathematics\\
Middle East Technical University, 
 Ankara-Turkey\\
	\texttt{yildiz.suleyman@metu.edu.tr} \\
     \And
     B\"ulent Karas\"ozen\\
     Institute of Applied Mathematics \& Department of Mathematics,
     Middle East Technical University\\
     Ankara-Turkey\\
     \texttt{bulent@metu.edu.tr}
}

\begin{document}
\maketitle

\date{Received: date / Accepted: date}

\begin{abstract}
Computationally efficient, structure-preserving reduced-order methods are developed for the Korteweg-de Vries (KdV) equations in Hamiltonian form. The semi-discretization in space by finite differences is based on the Hamiltonian structure. The resulting skew-gradient system of ordinary differential equations (ODEs) is integrated with the linearly implicit Kahan's method, which preserves the Hamiltonian approximately. We have shown, using proper orthogonal decomposition (POD), the Hamiltonian structure of the full-order model (FOM) is preserved by the reduced-order model (ROM). The reduced model has the same linear-quadratic structure as the FOM. The quadratic nonlinear terms of the KdV equations are evaluated efficiently by the use of tensorial framework, clearly separating the offline-online cost of the FOMs and ROMs. The accuracy of the reduced solutions, preservation of the conserved quantities, and computational speed-up gained by ROMs are demonstrated for the one-dimensional single and coupled KdV equations, and two-dimensional Zakharov-Kuznetsov equation with soliton solutions.
\end{abstract}


\section{Introduction}
\label{sec:intro}

Numerical integration of large scale dynamical systems is
computationally costly and requires a large amount of computer memory for
applications in real-time and many query solutions. The reduced-order methods (ROMs) have emerged as a powerful approach to reduce the computational effort by
constructing a low-dimensional linear subspace, that approximately
represents the solution to the high-dimensional system \cite{Benner17b,Quarteroni14,Hesthaven16}. Projection-based model reduction is one of the well-known and widely used ROM techniques, generally implemented using offline-online decomposition. Proper
orthogonal decomposition (POD) with Galerkin projection is one of the most standard methods to construct a reduced basis \cite{Berkooz93,Sirovich87}. During the offline stage, a set of reduced basis is extracted from a collection of high-fidelity solutions. In the online stage, the reduced solutions are computed in the reduced space, spanned by a set of basis functions that represents the main dynamics of the full-order model (FOM).

Many dynamical systems have some mathematical structures, such as symmetry, symplecticity, and energy preservation. Numerical integrators that inherit such properties are referred to as geometric numerical integrators or structure-preserving integrators \cite{Hairer16}. They produce stable and qualitatively better numerical solutions than standard general-purpose integrators. Various symplectic and  multisymplectic algorithms have been extended to Hamiltonian partial differential equations (PDEs) to preserve conservation laws.
When a Hamiltonian PDE is considered, the Galerkin projection-based POD-ROM is not able to preserve the desired physical quantities of the original system because the Hamiltonian structure of the original system may not be retained in the reduced dynamical system.  The reduced-order solutions may exhibit spurious and unphysical artifacts, leading to instabilities and qualitatively wrong solution behavior. Therefore, ROMs are preferred, that preserve the geometric structure and conserved quantities of FOMs. In the recent years, several structure-preserving reduced-order methods have been  developed for Lagrangian systems \cite{Carlberg13}, for port-Hamiltonian systems \cite{Chaturantabu16}, for dissipative Hamiltonian systems \cite{Afkham19}, for canonical \cite{Afkham17,Haasdonk19,Hesthaven20,Peng16,Karasozen18}, and for non-canonical Hamiltonian PDEs  \cite{Gong17,Miyatake19,Hesthaven18}.

In this paper, we develop an efficient structure-preserving ROMs for the Korteweg-de Vries (KdV) equation. The KdV equation is an integrable Hamiltonian PDE with a constant Poisson structure. The conserved quantities of the  KdV equation are the cubic Hamiltonian (energy), quadratic momentum and linear mass. The KdV equation is a nonlinear dispersive equation with smooth solutions.
There are relatively few papers concerning reduced-order modeling of the KdV equation. In \cite{Gerbau14} ROMs are constructed based on Lax-pairs, and in \cite{Hesthaven18} a greedy POD algorithm is developed with discrete empirical interpolation method (DEIM) based on the Poisson structure. In \cite{Miyatake19} structure-preserving POD and DEIM are constructed preserving first integrals of the KdV equation, and in \cite{Ehrlacher19} for one-dimensional conservative PDEs in Wasserstein space, ROMs are constructed including the KdV equation. For nonlinear PDEs without
polynomial structure, using hyper-reduction methods like the empirical interpolation (EIM) \cite{Maday14} and DEIM \cite{Chaturantabut10nmr}, the computational efficiency is discovered in solving the reduced system, i.e., in the online stage. When nonlinear PDEs like the KdV equation have polynomial structure, projecting the FOM onto the reduced space yields low-dimensional matrix operators that preserve the polynomial structure of the FOMs. Using the offline-online decomposition, computationally efficient ROMs can be constructed.

The KdV equation is discretized in space using various methods; finite difference, finite-volume, finite-element, spectral elements. Finite-volume and finite-element methods are suited for complex geometries, while spectral methods have higher order accuracy,  but lead to dense matrices for two-dimensional problems. Here, we consider only one-dimensional and rectangular two-dimensional domains in space. In this paper, we discretize the KdV equation in space by finite differences while preserving the skew-symmetry of the Poisson structure. The resulting skew-gradient system of ordinary differential equations (ODEs) preserves the energy, momentum, and mass at the discrete level. The resulting semi-discrete system is a linear-quadratic ODE system.
Most of the energy-preserving methods proposed so far are fully implicit methods, like the average vector field (AVF) method  \cite{Celledoni12}, where a system of nonlinear equations has to be solved at each time step by iterative methods like Newton's method or fixed-point iteration. The computational cost of the iterative solvers increases with the number of iterations and system size. The AVF method also requires the use of hyper-reduction techniques such as the DEIM to reduce the computational cost of the nonlinear terms in the ROMs  \cite{Karasozen18}.
For time discretization, we use as an alternative to AVF, the second-order linearly implicit Kahan's method \cite{Kahan97,Celledoni13} which is designed for ODEs with quadratic polynomial terms, obtained by semi-discretization of the KdV equation in space by finite differences.  In contrast to the fully implicit energy preserving schemes such as the average vector field  (AVF) method and the mid-point method, Kahan's method requires only one step Newton iteration at each time step for linear-quadratic systems such as the semi-discrete KdV equation  \cite{Celledoni13}. Kahan's method preserves the cubic integrals such as the Hamiltonians at the discrete-time level \cite{Celledoni15}. Applying POD in the tensorial framework (TPOD) \cite{Benner18,Benner15,Kramer19} by exploiting matricizations of tensors, the TPOD-ROM for the KdV equation with quadratic nonlinearity recovers an efficient offline-online decomposition.
The offline computation is accelerated by the use of tensor techniques  like matricizations of tensors \cite{Benner15,Benner18,Benner19,Kramer19}.
Here we make use of the sparse matrix technique MULTIPROD \cite{Leva08mmm} to further speed up the tensor calculations in the offline stage. We show the computational efficiency of the TPOD for three different KdV equations with soliton solutions; the one-dimensional single and coupled KdV equations, and the Zakharov-Kuznetsov equation which is a two-dimensional KdV equation.

The paper organized as follows. In Section \ref{fom} we introduce the FOM for three types of the KdV equations. In Section~\ref{rom} the structure-preserving ROMs with POD and TPOD are developed. We present in Section~\ref{num} numerical experiments demonstrating the preservation of the invariants accurately by ROMs with a low computational cost. The paper ends with concluding remarks in  Section~\ref{con}.
Through the paper, variables are denoted by plain letters, vectors are denoted by bold letters, and matrices and tensors are denoted by capital letters.

\section{Full-order model}
\label{fom}

KdV equation is a  dispersive, nonlinear hyperbolic equation with smooth solutions. It describes the propagation of long, one-dimensional waves,  including shallow-water waves, long internal waves in the ocean, ion-acoustic waves in a plasma, acoustic waves on a crystal lattice, and more. Dispersion and nonlinearity can interact to produce permanent and localized waveforms.
The KdV equation is a Hamiltonian PDE with a constant Poisson structure. It possesses bi-Hamiltonian structure \cite{Nutku90kdv,Karasozen13}, i.e., there exists an infinite number of invariants and therefore it is completely integrable. It was  solved using various geometric integrators; symplectic and multisymplectic methods \cite{Ascher05,Chen11,Reich01}, energy preserving integrators \cite{Karasozen13,Eidnes20,Karasozen12}.
In this section,  we construct FOMs by discretizing the one-dimensional single and coupled  KdV equations, and the two-dimensional KdV equation, i.e., Zakharov-Kuznetsov equation, in space and time.

\subsection{Single KdV equation}

The one-dimensional   KdV equation is given as
\begin{align} \label{kdv1}
\partial_t u &= -\alpha u\partial_xu - \mu \partial_{xxx}u ,
\end{align}
in a space-time domain $[a,b]\times [0,T]$ ($a<b, \; T>0$), with an initial condition and the periodic boundary condition
$$
u(x,0)= u^0(x), \qquad u(a,t) = u(b,t),
$$
with the real parameters $\alpha$ and $\mu$.
The KdV equation \eqref{kdv1} can be written as  a Hamiltonian PDE of the following form
\begin{equation*}
\partial_t u  = \mathcal{S} \frac{\delta \mathcal{H}}{\delta u},
\end{equation*}
where $\delta$ and $\partial$ denote the variational derivative and partial derivative, respectively.  The constant skew-adjoint operator (Poisson tensor) $\mathcal{S}$ and the Hamiltonian functional $\mathcal{H}$ are given by
\begin{equation*}
\mathcal{S} = \partial_x, \quad \mathcal{H}(u) = \int_a^b \left(-\frac{\alpha}{6}u^3 +  \frac{\mu}{2}(\partial_x u)^2   \right)dx.
\end{equation*}
The KdV equation \eqref{kdv1} is  completely integrable, i.e., it has infinitely many invariants. Among them, the momentum $\mathcal{I}_1= \int u^2 dx$, and the mass $\mathcal{I}_2= \int u dx$ are the most important ones.

Semi-discrete form of the KdV equation is obtained on the partition of the spatial interval $[a,b]$ into $N_x$ uniform elements
$$
a = x_1<x_2 < \cdots < x_{N_x} < x_{N_x+1}=b, \quad \Delta x = (b-a)/(N_x).
$$
Then we set semi-discrete solution vector as $\bm{u}:=\bm{u}(t)=(u_1(t),\ldots , u_{N_x}(t))^T$, where $u_i(t)= u(x_i,t)$, $i=1,\ldots ,N_x$.
The discrete Hamiltonian $H(\bm{u})$ is given by
\begin{equation}\label{dham1}
H(\bm{u}) = \sum_{i=1}^{N_x} \left(-\frac{\alpha}{6}u_i^3 + \frac{\mu}{2}\left( \frac{u_{i+1} -u_i}{\Delta x}\right)^2 \right)\Delta x.
\end{equation}
Similarly, the discrete momentum and mass are given as
\begin{equation*}
I_1(\bm{u}) = \sum_{i=1}^{N_x} u_i^2 \Delta x, \quad I_2(\bm{u}) = \sum_{i=1}^{N_x} u_i \Delta x.
\end{equation*}

The semi-discretized KdV equation \eqref{kdv1} is a Hamiltonian system of ODEs, equivalently  a skew-gradient system
\begin{equation} \label{hamode}
\bm{u}_t = S\nabla H(\bm{u}),
\end{equation}
with the discrete gradient $\nabla H(\bm{u})$ and the constant skew-symmetric matrix $S$
$$
\nabla H(\bm{u}) = -\frac{\alpha}{2}\bm{u}\odot\bm{u} - \mu D_2\bm{u}, \quad S=D_1,
$$
where $\odot$ denotes the element-wise multiplication of vectors.
The matrices $D_1\in\bbR^{N_x\times N_x}$ and $D_2\in\bbR^{N_x\times N_x}$ correspond to the centred  finite difference discretization of the first and second order derivative operators $\partial_x$ and $\partial_{xx}$, respectively, which are given under periodic boundary  conditions  by
\begin{equation} \label{difo}
 D_1:=\frac{1}{2\Delta x}
	\begin{pmatrix}
	0 & 1 & & & -1\\
	-1 & 0 & 1 & & \\
	  & \ddots &\ddots &\ddots & \\
	 &  & -1 & 0 & 1\\
	1 &  & & -1 & 0\\
	\end{pmatrix},
\;
D_2 := \frac{1}{\Delta x ^2}
	\begin{pmatrix}
-2 & 1 & & & 1\\
	1 & -2 & 1 & & \\
	  & \ddots &\ddots &\ddots & \\
	 &  & 1 & -2 & 1\\
	1 &  & & 1 & -2\\
	\end{pmatrix},
\end{equation}
where $D_1$ is skew-symmetric as an approximation of the skew-adjoint  Poisson tensor  ${\mathcal S}$.
Then, the semi-discretized KdV equation \eqref{kdv1} can be written as
\begin{equation} \label{kdv1sd}
	\bm{u}_t =  - \underbrace{ \mu D_3\bm{u}}_\text{linear} - \underbrace{\frac{\alpha}{2} D_1(\bm{u}\odot\bm{u})}_\text{quadratic},
\end{equation}
where the skew-symmetric matrix $D_3:=D_1D_2$ approximates the third order derivative $\partial_{xxx}$.

For time discretization, we divide the time interval $[0,T]$ into $N_t$ uniform elements $
0=t_0<t_1 < \cdots < t_{N_t}=T$, $\Delta t = T/N_t$, and we denote by $\bm{u}^k= \bm{u}(t_k)$ the full discrete approximation vector at time $t_k$, $k=0,\ldots , N_t$.
The semi-discrete KdV equation \eqref{kdv1sd} is a linear-quadratic  system of ODEs  of the following form
\begin{equation}\label{quad}
\bm{u}_t = \bm{f}(\bm{u}) := B_qQ(\bm{u})+B_l\bm{u},
\end{equation}
with the quadratic vector field  $Q(\bm{u}) = (\bm{u}\odot\bm{u})$  and the skew-symmetric matrices $B_l=-\mu D_3$ and $B_q=(-\alpha/2) D_1$.
As the time integrator, we use Kahan's method \cite{Celledoni12,Kahan97} whose application to the linear -quadratic system \eqref{quad} yields
\begin{equation*}
\frac{\bm{u}^{k+1} - \bm{u}^k}{\Delta t} = B_q\tilde{Q}(\bm{u}^k,\bm{u}^{k+1}) + \frac{1}{2}B_l(\bm{u}^k + \bm{u}^{k+1}) ,
\end{equation*}
where the symmetric bilinear form $\tilde{Q}(\cdot ,\cdot )$ is obtained by the polarization of the quadratic vector field $Q(\cdot)$ as follows \cite{Celledoni15}
\begin{equation*}
\tilde{Q}(\bm{u}^k,\bm{u}^{k+1}) := \frac{1}{2}\left( Q(\bm{u}^k+\bm{u}^{k+1}) - Q(\bm{u}^k) -  Q(\bm{u}^{k+1}) \right).
\end{equation*}

For a large class of Hamiltonian systems, the method
has a conserved quantity (related to energy) and an invariant \cite{Kahan97,SanzSerna94} .
Kahan's method is second order, time-reversal, and  linearly implicit for ODEs with quadratic vector fields \cite{Celledoni13} like the semi-discrete KdV equation
\eqref{quad},  i.e., $\bm{u}^{k+1}$ can be computed by solving a single linear system of equations
$$
\left ( I -\frac{\Delta t}{2} \bm{f}'(\bm{u}^k)\right ) \tilde{\bm{u}} = \Delta t \bm{f}(\bm{u}^k),\qquad \bm{u}^{k+1} = \bm{u}^k + \tilde{\bm{u}},
$$
where $I$ is the identity matrix and $\bm{f}'$ denotes the Jacobian matrix of $\bm{f}$.

Kahan's method is the restriction of a
Runge-Kutta method to quadratic vector fields
\cite{Celledoni13}
\begin{equation} \label{kahanrk}
\frac{\bm{u}^{k+1} - \bm{u}^k}{\Delta t} = -\frac{1}{2}\bm{f}(\bm{u}^k) + 2\bm{f} \left (\frac{\bm{u}^{k+1} + \bm{u}^k}{2}  \right )   - \frac{1}{2}\bm{f}(\bm{u}^{k+1}).
\end{equation}

Kahan's method preserves the Hamiltonian approximately, i.e.,  it preserves the modified Hamiltonian or the polarized energy
$$
\tilde{H}(\bm{u}) := H(\bm{u}) + \frac{1}{2}\Delta t \nabla H(\bm{u})^T (I-\frac{1}{2}\Delta t \bm{f}'(\bm{u}))^{-1} \bm{f}(\bm{u}),
$$
for all cubic Hamiltonian systems  with constant Poisson structure such as the KdV equation \cite{Celledoni13}.

Kahan's method has not been extensively studied for solving PDEs so far, with the
exception \cite{Kahan97}, where it is applied for solving the KdV equation.  It was shown that Kahan's  method exhibits  all favorable numerical properties like energy conservation, linear error growth with time. For Hamiltonian PDEs, using multiple points to
discretize the variational derivative, linearly implicit energy-preserving
schemes are defined \cite{Matsuo01}. These methods  are generalized
for deriving linearly implicit energy-preserving multistep methods for Hamiltonian
PDEs with polynomial invariants  \cite{Dahlby11}. A comparison of this approach and Kahan's
method applied to PDEs is given in \cite{Eidnes19}.
Recently a two-step generalization  of Kahan's method \cite{Eidnes20} is applied to  multisymplectic PDEs with cubic invariants. It was shown that discrete approximations to local
and global energy conservation laws are preserved  for  the one-dimensional KdV equation and the two-dimensional Zakharov-Kuznetsov equation.

Other energy preserving integrators like the implicit mid-point rule \cite{Miyatake19} and the AVF method \cite{Hesthaven18}, both are applied to the KdV equation in the context of reduced-order modelling, are fully implicit. The resulting nonlinear algebraic equations have to be solved by iteratively. We remark that implicit mid-point rule preserves only the quadratic Hamiltonians, whereas the AVF method preserves cubic Hamiltonians. For two-dimensional problems, where fully implicit schemes are computationally costly,  the
linearly implicit methods seem to provide for a competitive method.
The full order solutions can be speeded up in the periodic setting using the slit-step fast Fourier transformation (FFT) method which was originally proposed in \cite{Hardin73}.

\subsection{Coupled KdV equation}

As the second model, we consider the one-dimensional symmetric coupled KdV-KdV system \cite{Karasozen12,Bona07}
\begin{equation} \label{kdv2}
\begin{aligned}
\partial_t u  &= \frac{3}{2} u\partial_x u - \frac{1}{2} v\partial_xv -\partial_x v -\frac{1}{6}\partial_{xxx}v, \\
\partial_t v  &= -\partial_x u - \frac{1}{2}\partial_x (uv) -  \frac{1}{6}\partial_{xxx}u,
\end{aligned}
\end{equation}
which represents approximation to two-dimensional Euler equations for surface water waves propagation along a horizontal channel, where $u$ is the horizontal velocity and  $v$ is the deviation of the free surface from its rest position $x$. The initial and periodic boundary conditions are
$$
u_0(x,t) = u^0(x), \quad  v_0(x,t) = v^0(x), \qquad u(a,t)=u(b,t), \quad  v(a,t)=v(b,t).
$$
The corresponding Hamiltonian and skew-adjoint Poisson tensor  for the KdV-KdV system \eqref{kdv2} are given by
\begin{equation*}
\mathcal{H}(u,v) = \int_a^b  \left(-uv - \frac{1}{4}uv^2 - \frac{1}{4}u^3 - \frac{1}{6}u\partial_{xx}v \right)dx\; , \quad \mathcal{S} =
\begin{pmatrix}
\partial_x & 0\\
0 & \partial_x
\end{pmatrix}.
\end{equation*}
Additional invariants for the coupled KdV-KdV system \eqref{kdv2} are the momentum $\mathcal{I}_1= \int (u^2 +v^2)  dx$,  and the masses $\mathcal{I}_2= \int u dx$ and $\mathcal{I}_3= \int v dx$.
The discrete Hamiltonian $H(\bm{u},\bm{v})$ is given by
\begin{equation}\label{dham2}
H(\bm{u},\bm{v}) = \sum_{i=1}^{N_x} \left( -u_iv_i  - \frac{1}{4}u_iv_i^2  - \frac{1}{4}u_i^3 -  \frac{1}{6}u_i
\left( \frac{v_{i+1}-2v_i + v_{i-1}}{\Delta x^2} \right)  \right)\Delta x .
\end{equation}
The semi-discrete form of the coupled KdV-KdV system \eqref{kdv2} can be written as a skew-gradient system  with linear and quadratic terms
\begin{equation}\label{kdv2sd}
\begin{aligned}
	\bm{u}_t &= -\underbrace{ \left( D_1+\frac{1}{6}D_3 \right)\bm{v}}_\text{linear}
-\underbrace{\frac{3}{4} D_1 (\bm{u}\odot\bm{u}) - \frac{1}{4} D_1(\bm{v}\odot\bm{v})}_\text{quadratic} , \\
	\bm{v}_t &=  -\underbrace{\left( D_1+\frac{1}{6}D_3 \right)\bm{u}}_\text{linear} -  \underbrace{\frac{1}{2}D_1(\bm{u}\odot\bm{v})}_\text{quadratic}.
\end{aligned}
\end{equation}

\subsection{Zakharov-Kuznetsov equation}

The third model is the two-dimensional (2D) KdV equation known as the Zakharov-Kuznetsov equation  \cite{Iwasaki90,Nishiyama12,Zakharov74,Xu05}
\begin{equation} \label{zk}
\partial_t u  = -\alpha u\partial_x u -\mu(\partial_{xxx}u -\partial_{xyy} u),
\end{equation}
in the   space-time domain $([a,b]\times[c,d])\times [0,T]$ ($a<b,\; c<d,\; T>0$) with the initial condition and periodic boundary conditions
$$
u(x,y,0) = u^0(x,y), \quad  u(a,y,t) = u(b,y,t), \quad  u(x,c,t) = u(x,d,t).
$$
The skew-adjoint Poisson tensor and Hamiltonian are given as
	\begin{equation}\label{ham3}
	{\mathcal S} = \partial_x, \quad \mathcal{H}(u) =  \int_c^d\int_a^b  \left( -\frac{\alpha}{6}u^3  +  \frac{\mu}{2}\left ((\partial_x u)^2 + (\partial_y u)^2) \right) \right)  dxdy.
	\end{equation}
Additional invariants are  the momentum $I_1 = \iint \frac{1}{2}u^2dxdy$ and the mass $I_2= \iint u dxdy$. It describes the motion of nonlinear ion-acoustic waves in magnetized plasma.

For  space discretization, the spatial domain $\Omega =[a,b]\times [c,d]$ is divided into $N_x$ and $N_y$ elements in $x$ and $y$ directions, respectively, to form a rectangular mesh
\begin{align*}
&a =  x_1<x_2 < \cdots < x_{N_x} < x_{N_x+1}=b, &\Delta x = (b-a)/(N_x), \\
&c =  y_1<y_2 < \cdots < y_{N_y} < y_{N_y+1}=d, &\Delta y =  (d-c)/(N_y).
\end{align*}
Then, the semi-discrete solution vector is defined as
$$
\bm{u}:=\bm{u}(t)=(u_{1,1}(t),\ldots,u_{1,N_y}(t),u_{2,1}(t),\ldots ,  u_{N_x,N_y}(t))^T,
$$
where $u_{i,j}(t)= u(x_i,y_j,t)$, $i=1,\ldots , N_x$, $j=1,\ldots , N_y$.
The discrete form of the Hamiltonian in \eqref{ham3} is given by
\begin{equation} \label{ham3d}
H(\bm{u}) = \sum_{i=1}^{N_x}  \sum_{j=1}^{N_y} \left(-\frac{1}{6}(u_{i,j})^{3} +  \frac{\mu}{2}\left( \frac{u_{i+1,j} -u_{i,j}}{\Delta x}\right)^2  +
 \frac{\mu}{2}\left( \frac{u_{i,j+1} -u_{i,j}}{\Delta y}\right)^2 \right)\Delta x \Delta y.
\end{equation}
The semi-discrete form of the Zakharov-Kuznetsov equation \eqref{zk} is a skew-gradient system of the form
\begin{equation} \label{zksd}
\begin{aligned}
	\bm{u}_t &= S\nabla H(\bm{u}) = D_x\left( - \frac{\alpha}{2} (\bm{u}\odot\bm{u}) - \mu(D_{xx}+ D_{yy}) \bm{u} \right)\\
	         &= -\underbrace{\mu(D_{xxx}+ D_{xyy}) \bm{u}}_\text{linear} - \underbrace{\frac{\alpha}{2} D_x(\bm{u}\odot\bm{u})}_\text{quadratic} ,
\end{aligned}
\end{equation}
where we set $D_{xxx}:=D_xD_{xx}$, $D_{xyy}:=D_xD_{yy}$, and the 2D centred finite difference matrices $D_x,D_{xx},D_{yy}\in\bbR^{N_xN_y\times N_xN_y}$ are defined by
$$
D_x = D_1 \otimes I_y \; , \quad D_{xx} = D_2 \otimes I_y \; , \quad D_{yy} = I_x \otimes D_2,
$$
where $I_x$ and $I_y$ are $N_x$ and $N_y$ dimensional identity matrices, and the matrices $D_1$ and $D_2$ are the ones defined in \eqref{difo}, with appropriate dimension.

\section{Reduced-order model}
\label{rom}

Semi-discretization of KdV equations in Section~\ref{fom} leads to the following  system of linear-quadratic ODEs
\begin{equation} \label{quadfom}
\frac{d {\mathbf q}}{dt } = S \nabla_{\bm{q}} H(\bm{q}) = B_l{\mathbf q}  + B_qQ ({\mathbf q} ),
\end{equation}
where $\bm{q}\in \mathbb{R}^{N}$ is the state vector, $B_l,B_q \in \mathbb{R}^{N\times N}$  are the linear operators, $Q({\mathbf q} ):  \mathbb{R}^{N} \to \mathbb{R}^{N}$ is the quadratic operator, and $N$ is the  degree of freedom of the system , where $N=N_x$ for the single KDV system \eqref{kdv1sd}, $N=2N_x$ for the coupled KdV system \eqref{kdv2sd}, and $N=N_x\times N_y$ for the Zakharov-Kuznetsov system \eqref{zksd}.

The POD basis vectors are computed using the method of snapshots.
Consider the discrete state vector $\bm{q}$  as the solution to one of the KdV equations \eqref{kdv1sd}, \eqref{kdv2sd} or \eqref{zksd}. The snapshot matrix is defined as
$$
\mathcal{Q} := [ \bm{q}^1, \cdots , \bm{q}^{N_t} ] \in\mathbb{R}^{N\times N_t},
$$
where each column $\bm{q}^k\in\bbR^{N}$ is the full discrete solution vector at discrete time instances $t_k$, $k=1,\ldots ,N_t$. We then expand the singular value decomposition (SVD) of the snapshot matrix
$$
\mathcal{Q} = V\Sigma U^T,
$$
where the columns of $V\in\mathbb{R}^{N\times N_t}$ and $U\in\mathbb{R}^{N_t\times N_t}$ are the left and right singular vectors of $\mathcal{Q}$, respectively, and
$\Sigma \in\mathbb{R}^{N_t\times N_t}$ is the diagonal matrix whose diagonal elements are the singular values $\sigma_1 \ge \sigma_2 \ge \cdots\ge \sigma_{N_t}\ge 0$.

The $n$-POD basis matrix $V_n \in \mathbb{R}^{N\times n}$ minimizes the least squares error of the snapshot  reconstruction
\begin{equation*}
\min_{V_n \in \mathbb{R}^{N\times n}} ||\mathcal{Q}  -V_nV_n^T \mathcal{Q} ||^2_F=
\min_{V_n \in \mathbb{R}^{N\times n}}\sum_{k=1}^{N_t} ||\bm{q}^k  -V_nV_n^T \bm{q}^k ||^2_2 = \sum_{k=n+1}^{N_t} \sigma_k^2,
\end{equation*}
where $\|\cdot \|_2$  denotes the Euclidean $2$-norm and  $\|\cdot \|_F$  denotes the Frobenius norm. The optimal solution of basis matrix $V_n$ to this problem is given by
the $n$ left singular vectors of $\mathcal{Q}$  corresponding to the $n$ largest singular values.

The POD state approximation is $\bm{q} \approx \widehat{\bm{q}} = V_n  \bm{q}_r$, where $\bm{q}_r\in \mathbb{R}^n$ is the reduced state vector.  The POD reduced model is then defined by Galerkin
projection
\begin{align} \label{hamoder}
\frac{d}{d t} \bm{q}_r = V_n^T  S \nabla_{\bm{q}} H(V_n \bm{q}_r).
\end{align}
Although the matrix $S$ is a constant skew-symmetric matrix,
the reduced-order system \eqref{hamoder} based on Galerkin projection
  is not necessarily a skew-gradient system in general.
The Hamiltonian structure   can be preserved by inserting
$ V_nV_n^T\in\bbR^{N\times N}$ between $S$ and $\nabla_{\bm{q}} H(V_n  \bm{q}_r)$ in \eqref{hamoder}, which
yields a small skew-gradient system  \cite{Karasozen18,Gong17,Miyatake19}
\begin{equation} \label{kdvred}
\frac{\rmd}{\rmd t} \bm{q}_r = V_n^T S V_n V_n^T  \nabla_{\bm{q}} H(V_n\bm{q}_r)  = \widehat{S}\nabla_{\bm{q}_r} \widehat{H}(\bm{q}_r)
\end{equation}
where $\widehat{S} := V_n^T S V_n$ and $\widehat{H}(\bm{q}_r) := H(V_n \bm{q}_r) $.

The reduced cubic Hamiltonian $\widehat{H}$ is preserved by ROM,  because the ROM \eqref{kdvred} has the same skew-gradient form as the FOM \eqref{hamode}. We remark that the periodic boundary conditions in the FOM are preserved in the ROMs \cite{Sanderse20}.

The POD basis for the coupled PDEs, like the coupled KdV equation \eqref{kdv2sd} are usually computed by stacking all $\bm{u}$ and $\bm{v}$ in one vector $\bm{q}= (\bm{u},\bm{v})^T$ and by taking the SVD of the snapshot data. But the resulting ROMs do not preserve the coupling topology structure of the FOM \cite{Benner15,Reiss07,Kramer20} and  produce unstable reduced solutions. In order to maintain the coupling structure in ROMs, the POD basis vectors are computed separately for each the state vector $\bm{u}$ and $\bm{v}$. Let $ Q_u,Q_v\in \mathbb{R}^{N\times N_t} $ be snapshot matrices  for each state vector
$$
Q_u =  \left[\bm{u}^1 ,\ldots,\bm{u}^{N_t}\right], \quad  Q_v =  \left[\bm{v}^1 ,\ldots,\bm{v}^{N_t}\right].
$$
The POD basis are computed taking the SVD of the snapshot matrix $Q\in \mathbb{R}^{2N\times N_t}$
$$
Q =
\begin{pmatrix} Q_u\\Q_v
\end{pmatrix}
= \begin{pmatrix}V_u &  \\ &  V_v  \\    \end{pmatrix}
\begin{pmatrix} \Sigma_u &  \\  & \Sigma_v  \end{pmatrix}
\begin{pmatrix} U_u^T & \\  &  U_v^T  \end{pmatrix}.
$$

For PDEs like KdV equations with polynomial nonlinearities, ROMs do not require approximating the nonlinear terms through sampling   hyper-reduction methods. Reduced-order operators can be precomputed in the offline stage. Projection of FOM onto the reduced space yields low-dimensional matrix operators that preserve the polynomial structure of the FOM. This is an advantage because the offline-online computation is separated in contrast to the  hyper-reduction techniques like discrete empirical interpolation method, which may cause inaccuracies or instabilities in the ROM solutions in long term simulations. Recently, for PDEs with polynomial nonlinearities, the computationally efficient ROMs are constructed by the use of some tools from tensor theory and by matricizations of tensors \cite{Benner15a,Benner18,Benner19}.

The dimension of the ROM \eqref{kdvred} is supposed to be much smaller than  the dimension of the FOM \eqref{quadfom} ($n \ll N$) for an  efficient
online computation of the ROM.  But the computation of the quadratic terms of the reduced system still depends on the dimension of the FOM, with the computational cost of order $\mathcal{O}(nN)$ \cite{Stefanescu14}.
 This can be avoided by applying TPOD  and exploiting the   tensor  matricization. TPOD separates the full spatial variables from the reduced time variables, allowing fast  nonlinear term computations in  the online stage. Using the Kronecker product $\otimes$, the FOM \eqref{quadfom} can be written as the following linear-quadratic ODEs
\begin{equation} \label{quadfomK}
\frac{d {\mathbf q}}{dt } = S \nabla_{\bm{q}} H(\bm{q}) = B_l{\mathbf q}  + B_qW ({\mathbf q} \otimes {\mathbf q}),
\end{equation}
where $W\in\mathbb{R}^{N\times N^2}$ is the matricized tensor which satisfies the identity $W ( \bm{q} \otimes \bm{q} ) = \bm{q}\odot\bm{q}$.
The linear-quadratic structure of the FOM \eqref{quadfomK} is preserved by the ROM  \cite{Benner15a}
\begin{equation} \label{eq:kdv1rom}
\frac{d }{dt }\bm{q}_r = \widehat{B}_l\bm{q}_r  + \widehat{B}_q\widehat{W} (\bm{q}_r \otimes \bm{q}_r),
\end{equation}
where, for the single  KdV equation \eqref{kdv1}, $\widehat{B}_l$, $\widehat{B}_q$ and  $\widehat{W}$ are given  as
$$
\widehat{B}_l= -\mu\widehat{S}V_n^TD_2V_n, \quad  \widehat{B}_q = -\frac{\alpha}{2}\widehat{S},
$$
$$
\widehat{S}=V_n^TD_1V_n, \quad \widehat{W} = V_n^T W(V_n \otimes V_n).
$$
The ROMs of the coupled KdV equation \eqref{kdv2} and the Zakharov-Kuznetsov equation \eqref{zk} can be defined  similarly.

Using the TPOD, the computational cost of the reduced quadratic term in the ROM \eqref{eq:kdv1rom} becomes of order $\mathcal{O}(n^3)$ \cite{Stefanescu14}, i.e., the offline and online computations are separated.
On the other hand, TPOD requires the computation of the reduced tensor $\widehat{W}$ in the offline stage, but the explicit computation of $V_n \otimes V_n$ is inefficient because of the order ${\mathcal O}(n^2N^2)$ of the  computational complexity.
In order to avoid from this computational burden, $V_n \otimes V_n$ is computed in an efficient way using $W$ by $\mu$-mode matricizations of tensors \cite{Benner15}.
Recently algorithms are developed using tensor techniques to compute  $\widehat{W}$ by exploiting the particular structure of Kronecker product \cite{Benner18,Benner19}, wherein,  $\widehat{W}$ is computed without explicitly forming $W$ with the complexity of order $\mathcal{O}(n^3N)$ in contrast to the $\mu$-mode (matrix) computation. The reduced matrix $\widehat{W}$ can be given in MATLAB notation as follows	
	\begin{align}\label{goyal}
	\widehat{W} = V_n^T W(V_n \otimes V_n)
	=V_n^T
	\begin{pmatrix}
	V_n(1,:)\otimes V_n(1,:)\\
	\vdots\\
	V_n(N,:)\otimes V_n(N,:)
	\end{pmatrix},
	\end{align}
which utilizes the structure of $ W(V_n\otimes V_n) $, without  explicit construction of $W$. In \cite{Benner18,Benner19}
 the CUR matrix approximation  \cite{Mahoney09} of $ W(V_n\otimes V_n) $ is used to increase computational efficiency. Instead, here we make use of  the  "MULTIPROD" \cite{Leva08mmm} to increase the computational efficiency of  $\widehat{W}$ in the offline stage.
The MULTIPROD\footnote{\href{https://www.mathworks.com/matlabcentral/fileexchange/8773-multiple-matrix-multiplications-with-array-expansion-enabled}{https://www.mathworks.com/matlabcentral/fileexchange/8773-multiple-matrix-multiplications-with-array-expansion-enabled}} handles multiple multiplications of the multi-dimensional arrays via virtual array expansion. It is  a fast and memory efficient generalization for arrays of the MATLAB matrix multiplication operator. For any given two vectors $ \mathbf{a} $ and $ \mathbf{b} $, the Kronecker product satisfies
	\begin{equation*}
	(\text{vec}{(\mathbf{b}\mathbf{a}^\top)})^\top =(\mathbf{a}\otimes \mathbf{b})^\top
	=\mathbf{a}^\top\otimes \mathbf{b}^\top,
	\end{equation*}
where vec$ (\cdot) $ denotes the vectorization of a matrix. Using the above identity, the matrix $C= W(V_n\otimes V_n)\in\mathbb{R}^{N\times n^2}$ can be constructed as
	\begin{equation}\label{C}
	C(i,:)=(\text{vec}(V_n(i,:)^\top V_n(i,:))^\top , \ \ i\in\{1,2,\ldots,N\}.
	\end{equation}
Reshaping the matrix $ V_n\in \mathbb{R}^{N\times n} $ as $ \widetilde{V}_n \in \mathbb{R}^{N\times 1 \times n} $ and computing MULTIPROD of $ V_n $ and $\widetilde{V}_n $ in the $2$nd and $3$rd dimensions, we obtain that
$$
\mathcal{C} =\text{MULTIPROD}( V_n,\widetilde{V}_n)\in \mathbb{R}^{N\times n \times n},
$$
where the matrix $C$ is recovered by reshaping the 3-dimensional array $\mathcal{C}$ into a matrix of dimension $N\times n^2$. Without MULTIPROD, the computation of the matrix $C$ in \eqref{C} requires $N$ for loops within each iteration the matrix product of two matrices of sizes $n\times 1$ and $1\times n$ are done. But, with the MULTIPROD, the matrix products are computed simultaneously in a single loop, and the matrix $\widehat{W}$ in \eqref{goyal} can be efficiently computed \cite{karasozen21}.

\section{Numerical results}
\label{num}

In this section, we demonstrate the performance of the structure-preserving ROM for the single KdV equation \eqref{kdv1} with one and two solitons, the coupled symmetric KdV-KdV  system \eqref{kdv2}, and  the Zakharov-Kuznetsov equation \eqref{zk}. For all the problems, we prescribe periodic boundary conditions on the given spatial domain.
In numerical test examples, we show only the preservation of the cubic integrals like the Hamiltonian (energy). Momentum as a quadratic invariant is preserved by all the Runge Kutta methods of type \eqref{kahanrk} including the Kahan's method and the implicit-midpoint rule. Linear invariants like the mass are automatically preserved by the Runge-Kutta methods.

All the simulations are performed on a machine with Intel CoreTM i7 2.5 GHz 64 bit CPU, 16 GB RAM, Windows 10, using
64 bit MatLab R2014.
 The snapshot matrices resulting from the space-time discretization of the KdV equations are large, making SVD computations costly. Therefore,  we use the randomized SVD (rSVD) algorithm  \cite{Halko11a} that  performs SVD of small matrices, to efficiently
generate a reduced basis.

In all examples, the  number of (POD) modes is determined by the relative information content (RIC) formula
\begin{equation} \label{ric}
\mathcal{E}_{\text{ric}}(n)=\left( \frac{\sum_{k=1}^{n} \sigma_{k}^2}{\sum_{k=1}^{N_t} \sigma_{k}^2  }\right) \times 100,
\end{equation}
which can be thought as the percentage energy captured from the FOM. According to the RIC   formula \eqref{ric}, we set the number of POD modes as the smallest positive integer $n$ satisfying $\mathcal{E}_{\text{ric}}(n)\geq 99.99$.

The accuracy of the ROM solutions are measured by the time averaged relative $L_2$-errors
\begin{equation}\label{relerr}
\|\bm{q}-\widehat{\bm{q}}\|_{\text{rel}}  =  \frac{1}{N_t}\sum_{k=1}^{N_t}\frac{\|\bm{q}^k-\widehat{\bm{q}}^k\|_{L^2(\Omega)}}{\|\bm{q}^k\|_{L^2(\Omega)}}, \quad
\|\bm{q}^k\|_{L^2(\Omega)}^2  =  \sum_{i=1}^N(\bm{q}_{i}^k)^2\Delta x\Delta y.
\end{equation}
We measure the preservation of the reduced conserved quantities using the time-averaged absolute  errors between the full and reduced quantities
\begin{equation} \label{conserr1}
\|E-\widehat{E}\|_{\text{abs}} = \frac{1}{N_t}\sum_{k=1}^{N_t}|E(\bm{q}^k)-\widehat{E}(\bm{q}_r^k)|,  \quad E= H,I_1,
\end{equation}
where $\widehat{E}(\bm{q}_r^k)=E(V_n\bm{q}_r^k)$ denotes the reduced quantity at the time $t_k$.

\subsection{Single KdV equation}

We consider the one-dimensional single KdV equation \eqref{kdv1} with $\alpha =6, \; \mu =1$  in the space-time domain $[-10,10] \times [0,50]$. For a positive parameter $\beta$, the initial condition  is set to $u(x,0) = \beta \sech^2 (\sqrt{\beta}x/2)$, which leads to one soliton solutions.
We set mesh size in space as $\Delta x =0.002$ and time step size is $\Delta t = 0.005$. The size  of the snapshot matrix is $Q\in \mathbb{R}^{10000\times 10000}$.

\begin{figure}[H]
\centering
\includegraphics[width=0.55\linewidth]{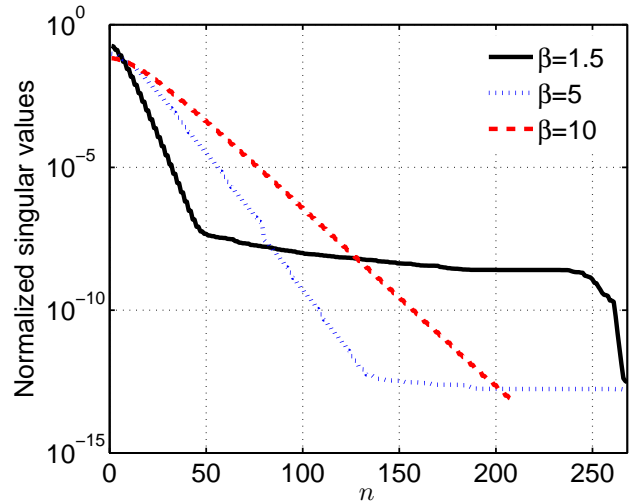}
\caption{Singular values of the snapshot matrices by different values of $\beta$.}
\label{ex1sing}
\end{figure}

The singular values decay much slowly for larger values of $\beta$  in Figure~\ref{ex1sing}. Consequently, more modes are needed  for  accurate computation of the reduced solutions with increasing  $\beta$. This behavior is characteristic for PDEs like the KdV equation exhibiting wave propagation phenomena, which require sufficiently large reduced spaces \cite{Ohlberger16}.

\begin{figure}[H]
\centering
\includegraphics[width=0.45\columnwidth]{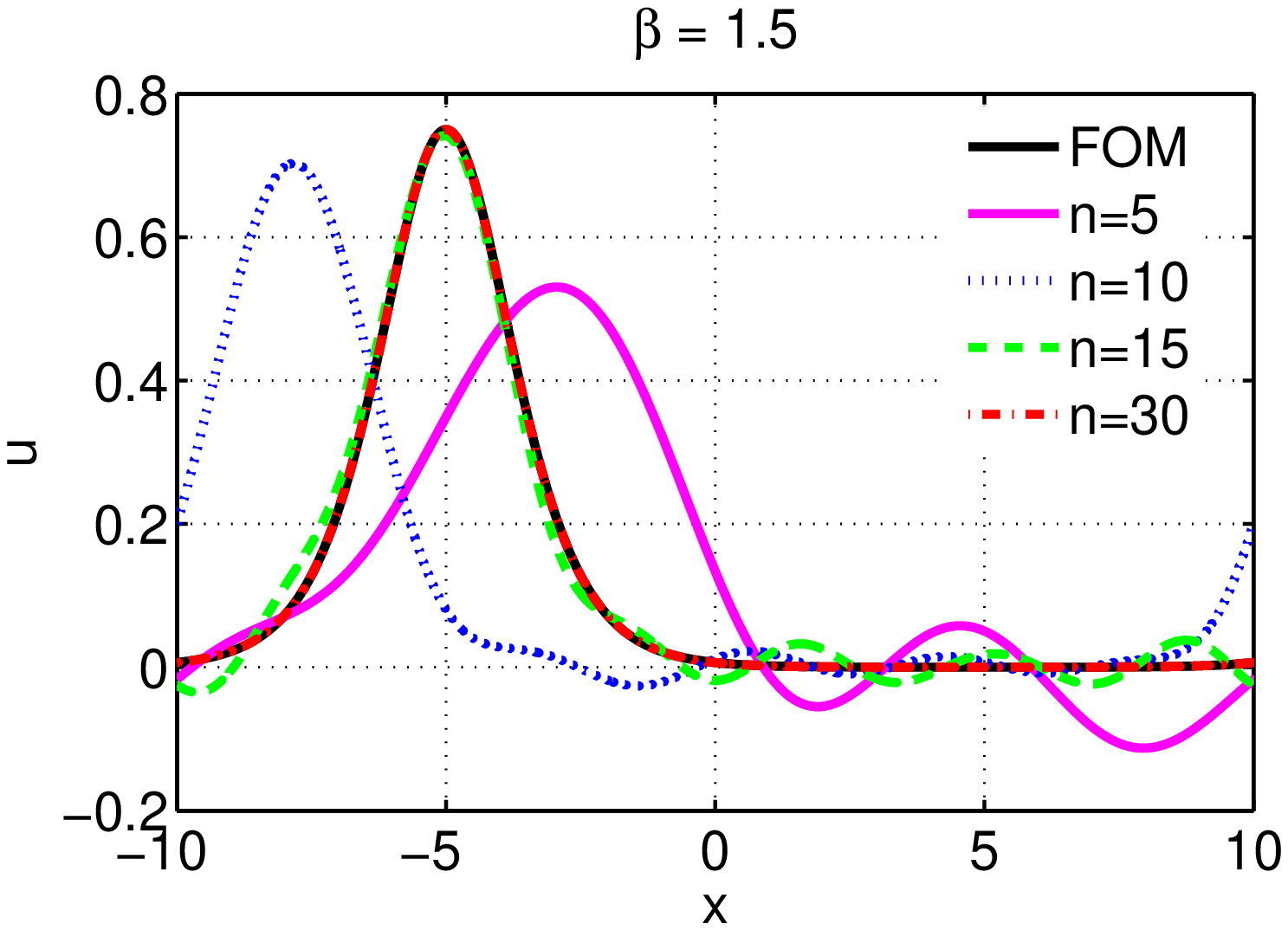}
\includegraphics[width=0.45\columnwidth]{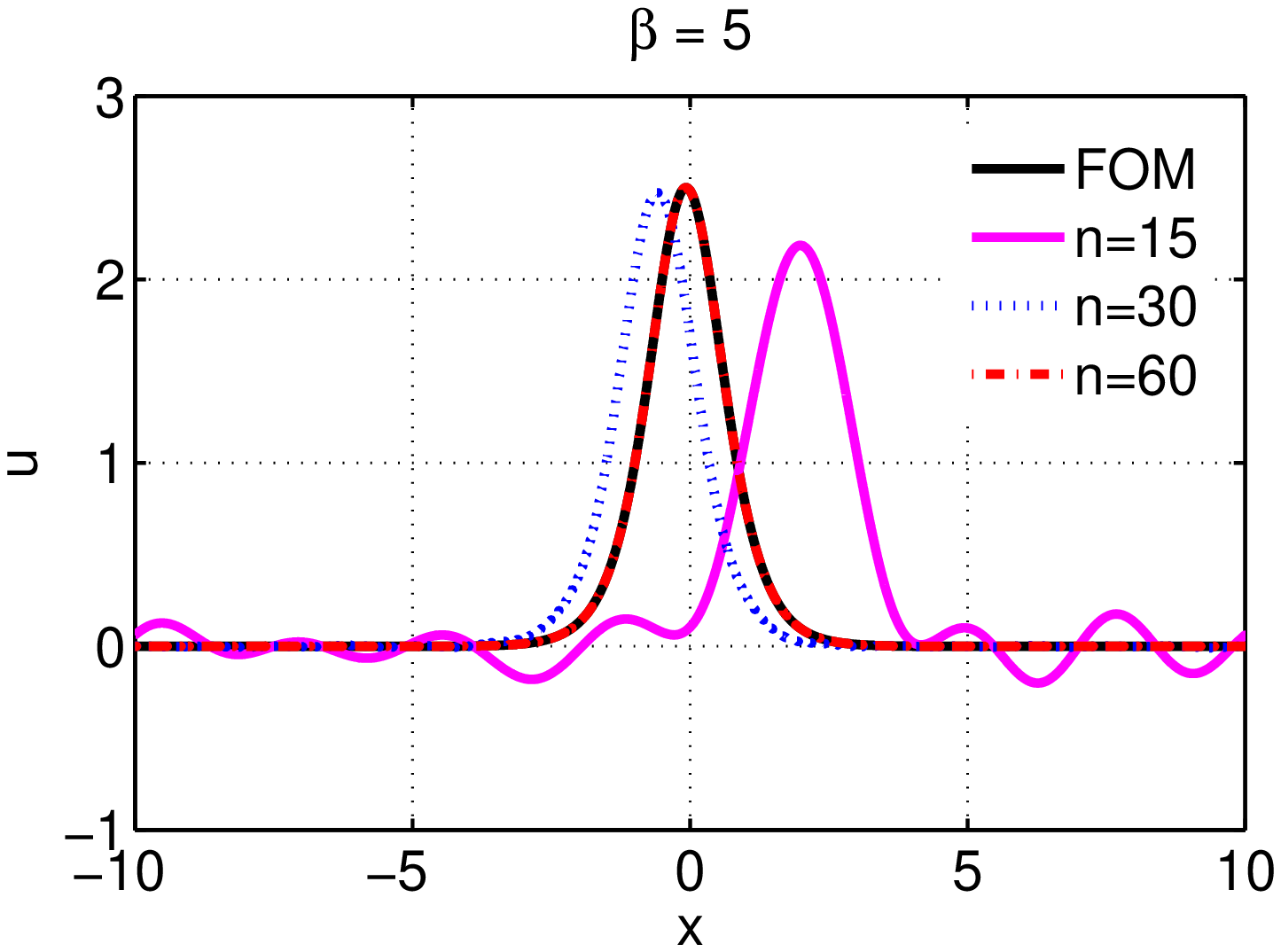}

\includegraphics[width=0.45\columnwidth]{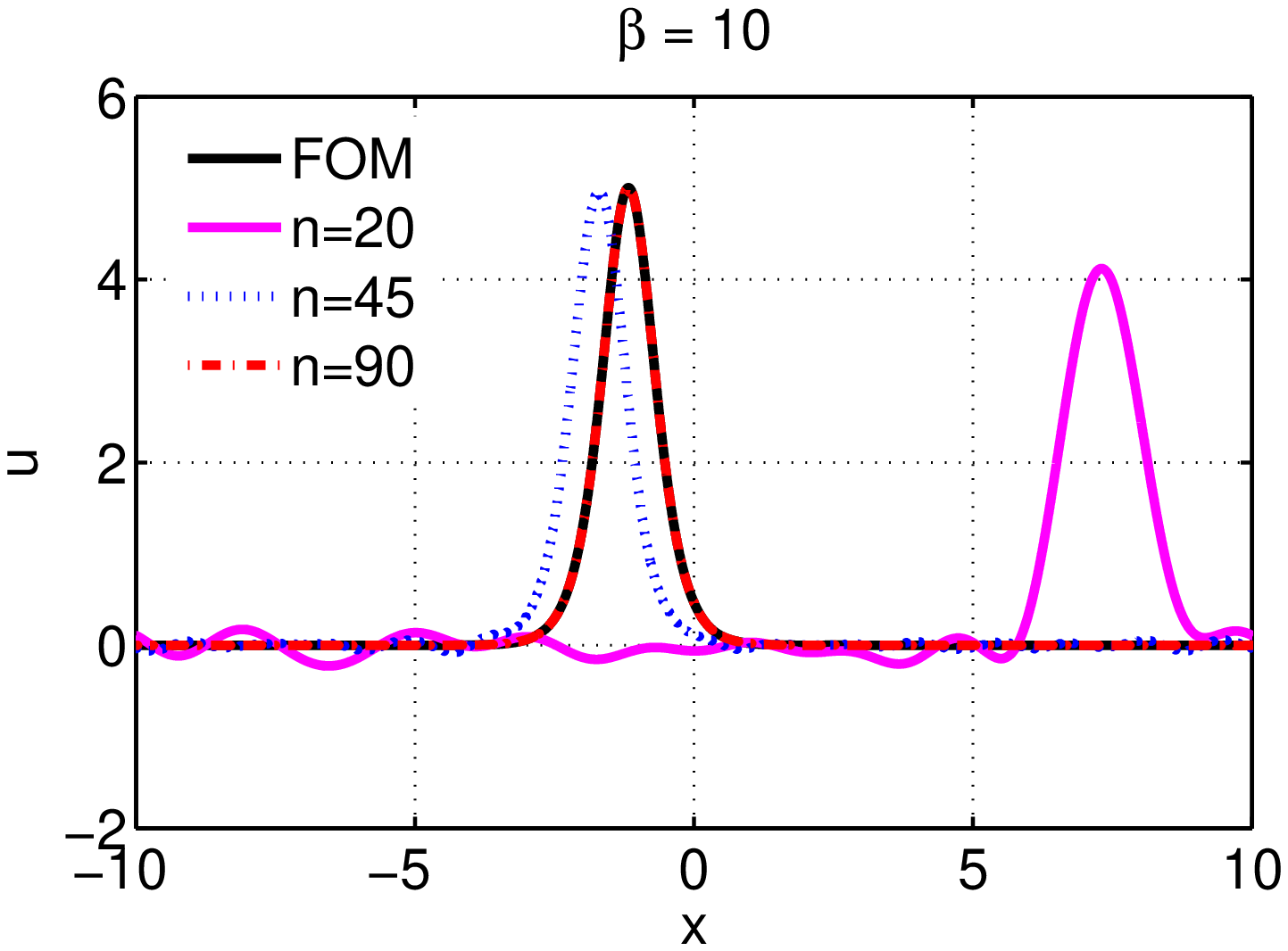}
\includegraphics[width=0.45\columnwidth]{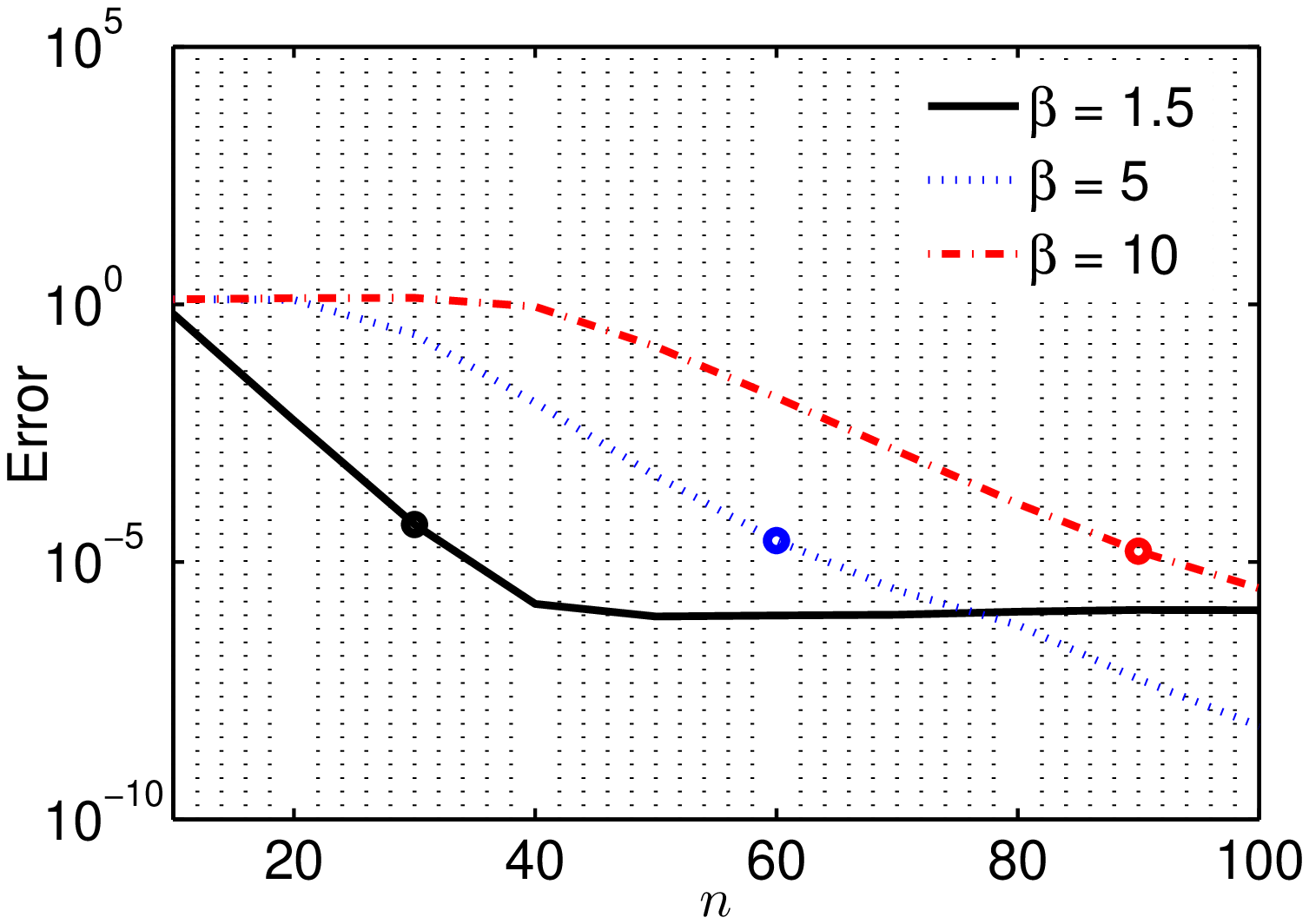}
\caption{ROM profiles at $T=50$ and relative solution errors \eqref{relerr} between FOMs and ROMs for different number of modes $n$ and for different values of $\beta$. The circles in the  bottom-right plot indicate the number of modes calculated according to the RIC formula \eqref{ric}.}
\label{ex1romsol}
\end{figure}

According to the RIC formula \eqref{ric}, the number of modes are taken as $n=30,60,90$ for  $\beta=1.5,5,10$, respectively.  In Figure \ref{ex1romsol} the reduced approximations are plotted for $\beta=1.5,5,10$ for increasing number of modes.  We observe that the relative $L^2$-errors  \eqref{relerr} between the full  and the reduced solutions  decrease as the number of modes increases in Figure \ref{ex1romsol}, bottom-right.
The accuracy of the reduced solutions is improved as the
number of modes is increased, upper and bottom left plots in Figure~\ref{ex1romsol}.
They are visually not distinguishable\footnote{Animations are available as the supplementary material "Ex1\_sol.mp4".} from the full solutions
for the number of modes selected by the RIC formula and indicated
by a circle in Figure~\ref{ex1romsol}, bottom-right.

\begin{figure}[H]
\centering
\includegraphics[width=0.9\columnwidth]{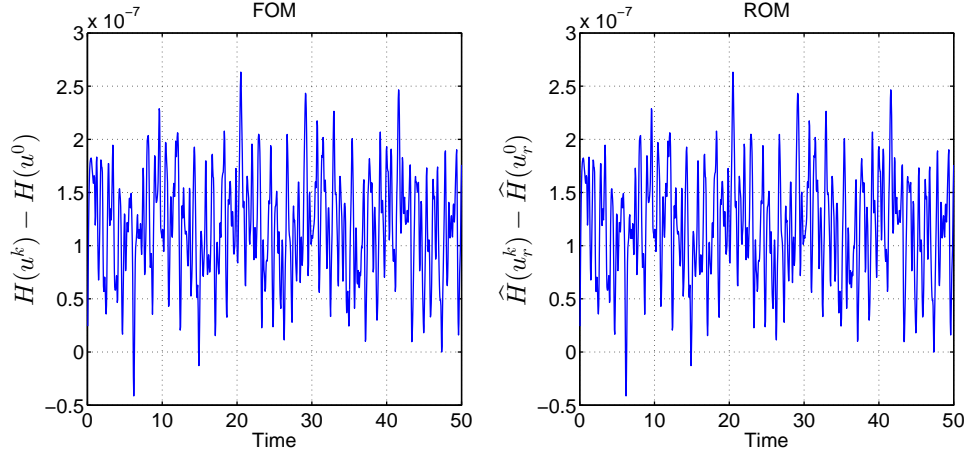}
\caption{Time evolution of the full (left) and the reduced (right) Hamiltonian errors.}
\label{ex1ham}
\end{figure}

\begin{figure}[H]
\centering
\includegraphics[width=0.55\columnwidth]{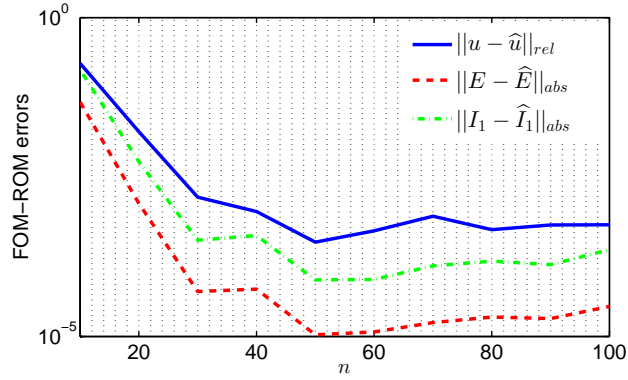}
\caption{Relative solution errors and absolute  errors  of the conserved quantities.  }
\label{ex1errors}
\end{figure}

Figure \ref{ex1ham} shows that the discrete cubic Hamiltonian \eqref{dham1} is preserved by the ROMs with high accuracy over time. The structure-preserving feature of the ROMs is well demonstrated by the solution errors \eqref{relerr} and errors of the conserved quantities \eqref{conserr1} in Figure~\eqref{ex1errors}. The relative FOM-ROM errors of the solutions and the errors in the Hamiltonian $H$ and the momentum $I_1$ are decreasing for an increasing number of modes with small oscillations around $n=50-80$.

In Table~\ref{kdv1beta} the relative solution errors \eqref{relerr},  conservation errors \eqref{conserr1} of the Hamiltonian and the momentum are given for $\beta= 1.5,5,10$.  With increasing values of $\beta$, more modes are needed for accurate reduced solutions and for the conservation of the Hamiltonian and the momentum.

\begin{table}[H]
	\caption{Hamiltonian, momentum and solution errors between the FOMs and ROMs \label{kdv1beta}}
\smallskip
		\resizebox{\textwidth}{!}{\begin{tabular}{|l|l|l|l|l|l|l|l|l|l|}  \hline
				& \multicolumn{3}{c|}{$\beta=1.5$}  & \multicolumn{3}{c|}{$\beta=5$}    & \multicolumn{3}{c|}{$\beta=10$}   \\ \hline
				\# modes  & $\|\bm{u}-\widehat{\bm{u}}\|_{\text{rel}}$ & $\|H-\widehat{H}\|_{\text{abs}}$ & $\|I_1-\widehat{I}_1\|_{\text{abs}}$ & $\|\bm{u}-\widehat{\bm{u}}\|_{\text{rel}}$ & $\|H-\widehat{H}\|_{\text{abs}}$ & $\|I_1-\widehat{I}_1\|_{\text{abs}}$ & $\|\bm{u}-\widehat{\bm{u}}\|_{\text{rel}}$ & $\|H-\widehat{H}\|_{\text{abs}}$ & $\|I_1-\widehat{I}_1\|_{\text{abs}}$ \\ \hline
				10  & 6.52e-01  & 1.34e-02  & 5.06e-03  & 1.26e+00  & 2.76e+00  & 6.06e-01  & 1.24e+00  & 2.98e+01  & 4.35e+00  \\ \hline
				20  & 5.53e-03  & 3.42e-05  & 4.84e-06  & 1.22e+00  & 2.04e-01  & 2.12e-02  & 1.31e+00  & 6.37e+00  & 4.90e-01  \\ \hline
				30  & 5.28e-05  & 4.43e-08  & 2.96e-09  & 2.60e-01  & 8.74e-03  & 4.69e-04  & 1.33e+00  & 8.51e-01  & 3.55e-02  \\ \hline
				40  & 1.52e-06  & 7.16e-11  & 1.04e-10  & 1.23e-02  & 2.79e-04  & 7.54e-06  & 8.94e-01  & 1.03e-01  & 1.45e-03  \\ \hline
				50  & 8.68e-07  & 1.02e-10  & 1.07e-10  & 4.64e-04  & 7.33e-06  & 6.00e-08  & 1.51e-01  & 8.10e-03  & 6.39e-05  \\ \hline
				60  & 9.13e-07  & 9.71e-11  & 1.07e-10  & 2.57e-05  & 1.69e-07  & 1.29e-09  & 1.56e-02  & 6.05e-04  & 2.09e-05  \\ \hline
				70  & 9.44e-07  & 9.42e-11  & 1.08e-10  & 2.83e-06  & 3.53e-09  & 8.49e-11  & 1.42e-03  & 4.09e-05  & 2.88e-06  \\ \hline
				80  & 1.08e-06  & 6.73e-11  & 1.09e-10  & 5.90e-07  & 2.41e-10  & 2.25e-12  & 1.35e-04  & 2.38e-06  & 3.39e-07  \\ \hline
				90  & 1.17e-06  & 5.86e-11  & 1.09e-10  & 5.33e-08  & 6.41e-12  & 3.23e-12  & 1.61e-05  & 1.53e-07  & 2.66e-08  \\ \hline
				100 & 1.15e-06  & 5.84e-11  & 1.09e-10  & 6.42e-09  & 7.75e-12  & 3.12e-12  & 3.06e-06  & 8.86e-09  & 2.18e-09  \\ \hline
		\end{tabular}}
\end{table}

\subsection{Two soliton interaction}

As the second test problem, we consider for $\alpha =1$ and $\mu =1$ the one-dimensional  two soliton KdV equation  \eqref{kdv1} with the
exact solution \cite{Brugnano19,Cai20,Liu16}
\begin{equation}\label{soliton}
u_e(x,t) = 12\frac{k_1^2e^{\xi_1} + k_2^2e^{\xi_2} +2(k_2-k_1)e^{\xi_1 +\xi_2} +\rho^2 (k_2^2e^{\xi_1} + k_1^2e^{\xi_2}e^{\xi_1 +\xi_2}}
{(1 + e^{\xi_1} +e^{\xi_2}+ \rho^2e^{\xi_1 +\xi_2})^2  }.
\end{equation}
The parameters are
\begin{eqnarray*}
k_1 & = & 0.4, \quad k_2 = 0.6, \quad \rho = (k_1-k_2)/(k_1+k_2) = -0.2,\\
\xi_1 & = & k_1x -k_1^3t + 4, \quad  \xi_1  =  k_2x -k_2^3t + 15.
\end{eqnarray*}
We take the space domain $\Omega = [-40, 40]$, and set the final time  $T=120$ as in \cite{Cai20,Liu16}.

To determine the experimental orders of convergence (EOC) of the high-fidelity solutions, the mesh size is uniformly refined  by a factor of two in both space and time dimensions. The EOC is calculated as

\begin{equation}\label{orderdef}
\mbox{order}=\frac{1}{\mbox{log} 2}\mbox{log}\left( \frac{\mbox{error}(\Delta x,\Delta t )}{\mbox{error}(\Delta x/2,\Delta t/2)}\right),
\end{equation}
where $\mbox{error}(\Delta x,\Delta t)$ denotes the relative $L^2$-error between the exact solution \eqref{soliton} and the numerical solution at the final time, computed with the spatial and temporal mesh sizes $\Delta x$ and $\Delta t$, respectively. The calculated
errors and their EOC are  summarized in Table~\ref{order}. They
confirm the expected second order rate of convergence of the centred finite difference scheme and Kahan's method.

\begin{table}[H]
\caption{Relative $L^2$-errors between the exact and FOM solutions, and  experimental order of convergence \label{order}}
\centering
\begin{tabular}{llllllll} \hline
$\Delta t $ & 0.5 & 0.25 & 0.125 & 0.0625 & 0.03125 & 0.016625 \\ \hline
$\Delta x $  & 4 & 2 & 1 & 0.5 & 0.25 & 0.125   \\ \hline
\hbox{Error} & 2.42e-00  &  9.68e-01 & 2.35e-01  &   5.72-02 & 1.42e-02 & 3.55e-03    \\ \hline
\hbox{Order} & - &  1.3226 &  2.0445  &  2.0359 & 2.0124 &   1.9970 \\ \hline
\end{tabular}
\end{table}

\begin{figure}[H]
\centering
\includegraphics[width=0.5\linewidth]{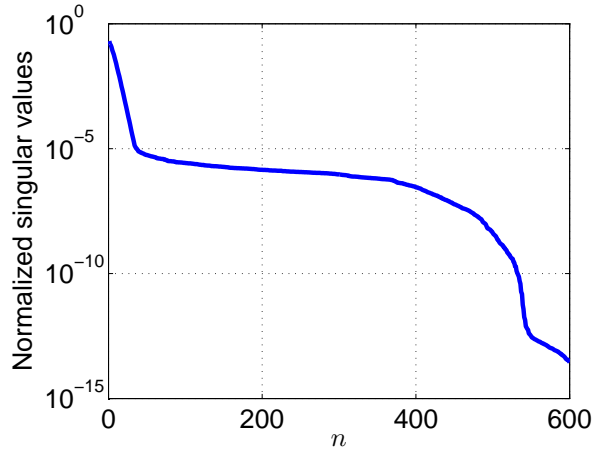}
\caption{Singular values of the snapshot matrix.}
\label{ex1singsoliton}
\end{figure}

We take the spatial mesh size as $\Delta x=0.125$ and  the time step $\Delta t=0.05$, that leads to the snapshot matrix $Q\in \mathbb{R}^{640\times 2400}$. The singular value spectrum in Figure \ref{ex1singsoliton} behaves similar to the single KdV equation with $\beta=1.5$. $30$  modes are sufficient to capture the behavior of the FOM soliton waves according to the RIC formula \eqref{ric}.

\begin{figure}[H]
\centering
\includegraphics[width=0.49\columnwidth]{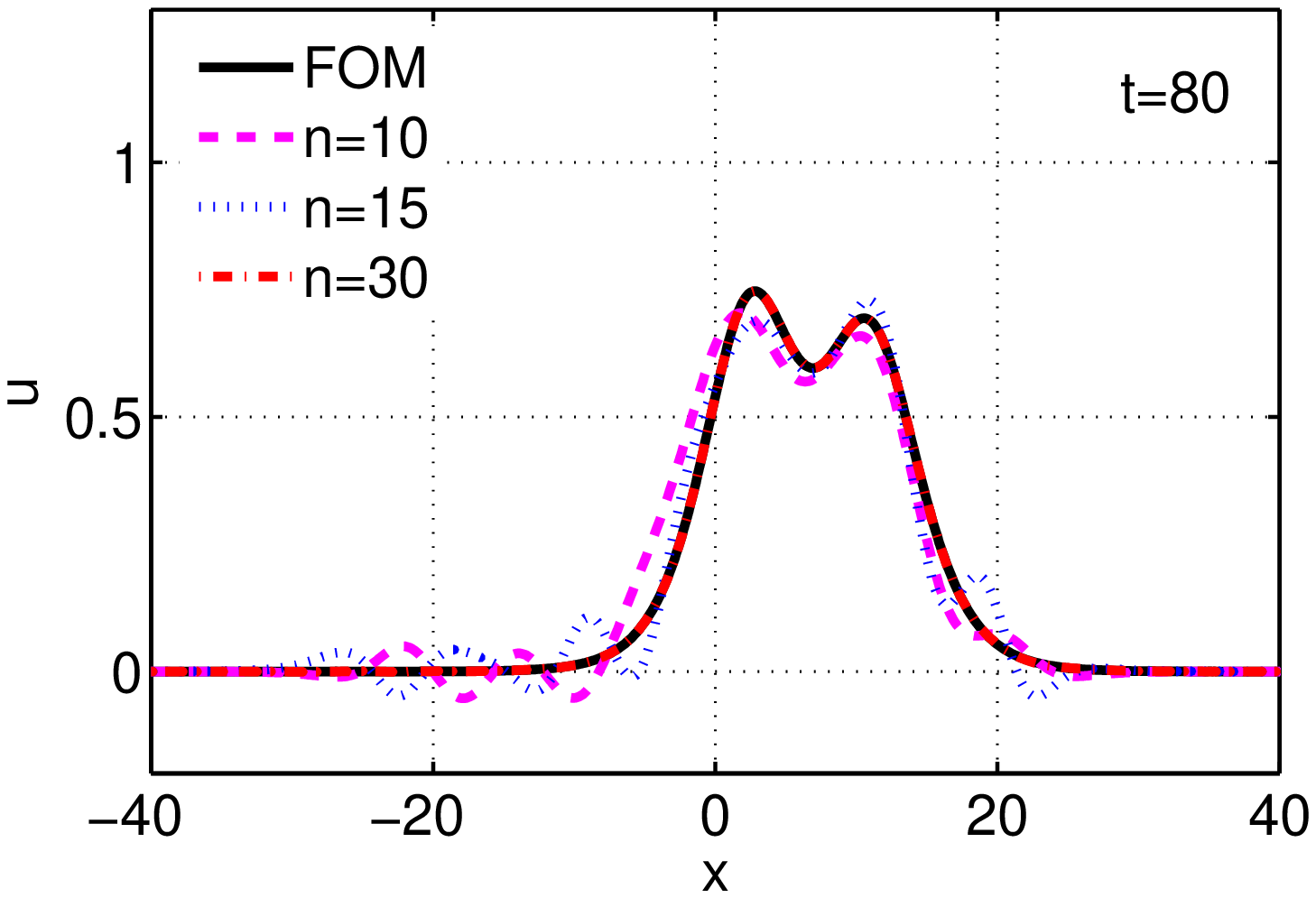}
\includegraphics[width=0.49\columnwidth]{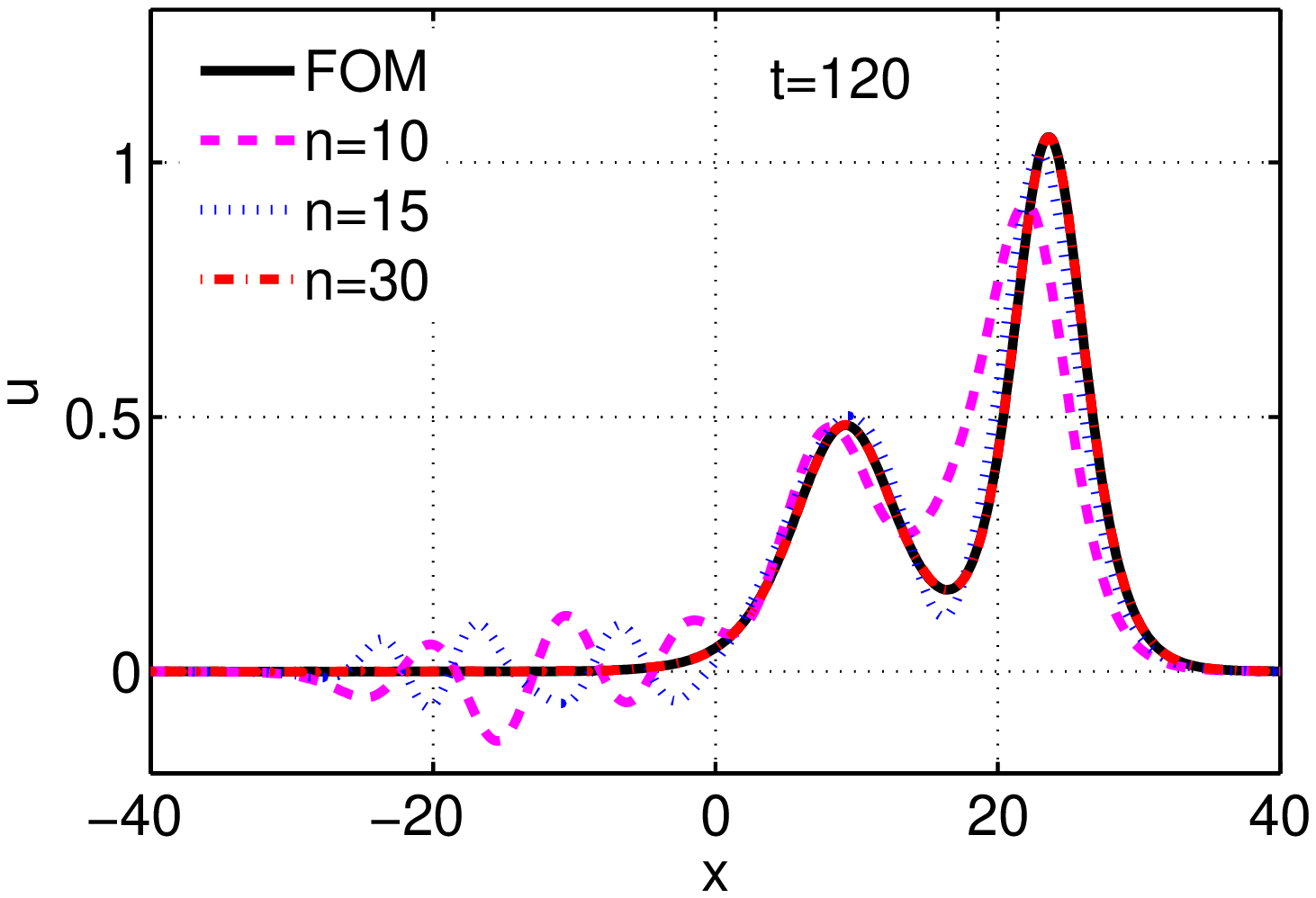}
\caption{FOM and ROM profiles at collision time $t=80$ and at final time $t=120$ for  different number of modes.}
\label{ex2_solitons}
\end{figure}

In Figure~\ref{ex2_solitons}, the two soliton waves\footnote{Animations are available as the supplementary material "Ex2\_sol.mp4".} with a taller and a lower one,  moving to the  right and  collide at $ t = 80$, continue moving away from each other until the final time $t= 120$ in Figure~\ref{ex2_solitons} as in \cite{Cai20,Liu16}. The ROM profiles  in Figure~\ref{ex2_solitons} at the collision time and at the final time show that with an increasing number of modes  they approximate the full solutions more closely  and  finally catch them  for $n=30$ according to the RIC formula \eqref{ric}.

\begin{figure}[H]
\centering
\includegraphics[width=0.6\columnwidth]{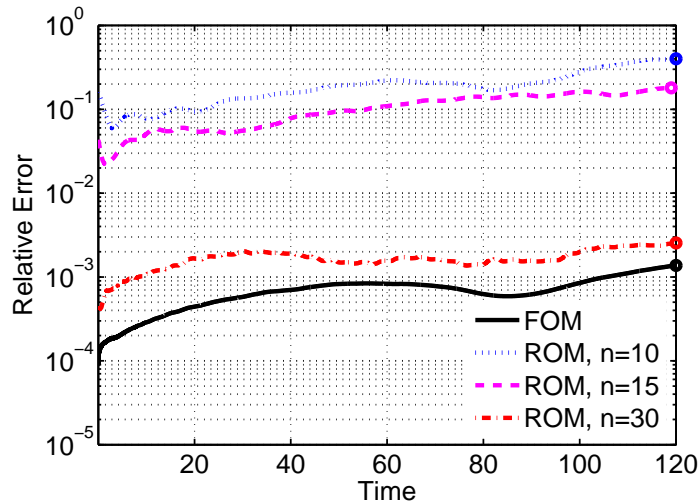}
\caption{Error propagation of relative $L^2$-errors between exact and FOM/ROM solutions.}
\label{ex2_maxerror}
\end{figure}

Furthermore, we show the propagation of the relative $L^2$-errors $\|u(x,t)-u_e(x,t)\|_{L^2(\Omega)}/\|u_e(x,t)\|_{L^2(\Omega)}$ between the exact solution
$u_e(x,t)$ and FOM/ROM solutions in Figure~\ref{ex2_maxerror}. The circles indicate that the maximum of the errors occur at the final time. In the reduced order modelling framework, the reduced solutions  are expected to behave similar to the full solutions,  since the reduced space is constructed from the FOM. Correspondingly, the errors of the reduced solutions in Figure~\ref{ex2_maxerror} show similar behavior  as the full  solution errors. This also indicates that the location and the shape of the full soliton waves are well-captured by the ROM solutions with increasing number of modes. The almost linear error growth rate in time in  Figure~\ref{ex2_maxerror} is characteristic for Hamiltonian preserving and for the  geometric integrators  \cite{Hairer06} including Kahan's method.

\begin{figure}[H]
\centering
\includegraphics[width=0.95\columnwidth]{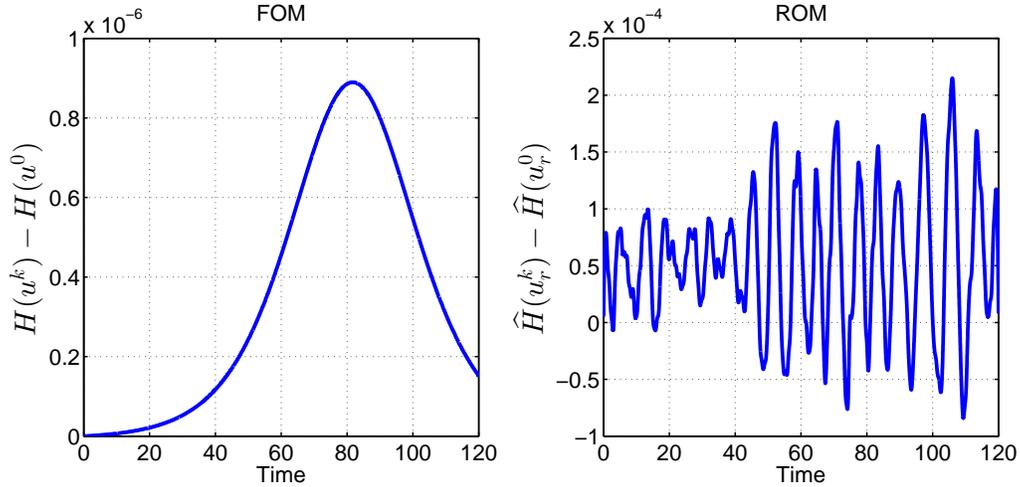}
\caption{Time evolution of the full (left) and the reduced (right) Hamiltonian errors.}
\label{ex1_energysoliton}
\end{figure}

In Figure~\ref{ex1_energysoliton}, the Hamiltonian errors do not show any drift, they are
preserved not with a high accuracy as  in  the single soliton example, which might be due to the interaction of the solitons.

\subsection{Coupled KdV equation}

Symmetric  KdV-KdV equation under periodic boundary conditions possesses solitary pulse solutions decaying symmetrically to oscillations of
small, constant amplitude \cite{Bona07,Bona08}.  The solutions are in the form of
traveling waves with main pulses like the classical solitary waves and dispersive oscillations
following the main pulses. For the coupled KdV-KdV equation \eqref{kdv2}, we take the initial conditions as in  \cite{Karasozen12,Bona08}
$$
u(x,0)=0 \; , \quad v(x,0)= 0.3e^{-(x+100)^2/25}.
$$
We set the space-time  domain as $[-150,150]\times [0,50]$, and  the mesh sizes are $\Delta x=0.1$ and $\Delta t=0.05$. The size of the snapshot matrix is
$Q\in \mathbb{R}^{3000\times 1000}$.

\begin{figure}[H]
\centering
\includegraphics[width=0.5\linewidth]{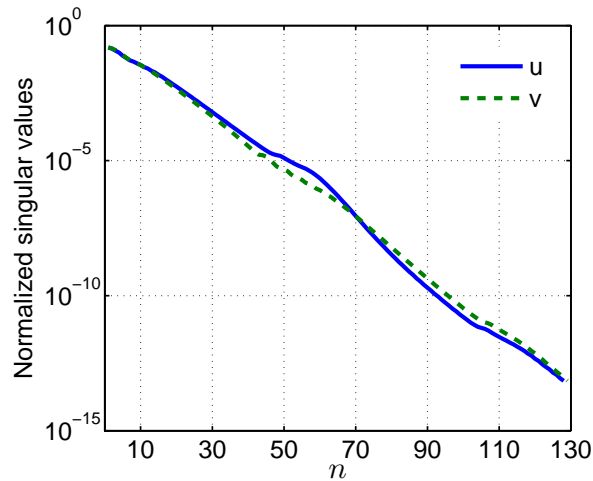}
\caption{Singular values of the snapshot matrices.}
\label{ex2sing}
\end{figure}

\begin{figure}[H]
\centering
\includegraphics[width=0.95\linewidth]{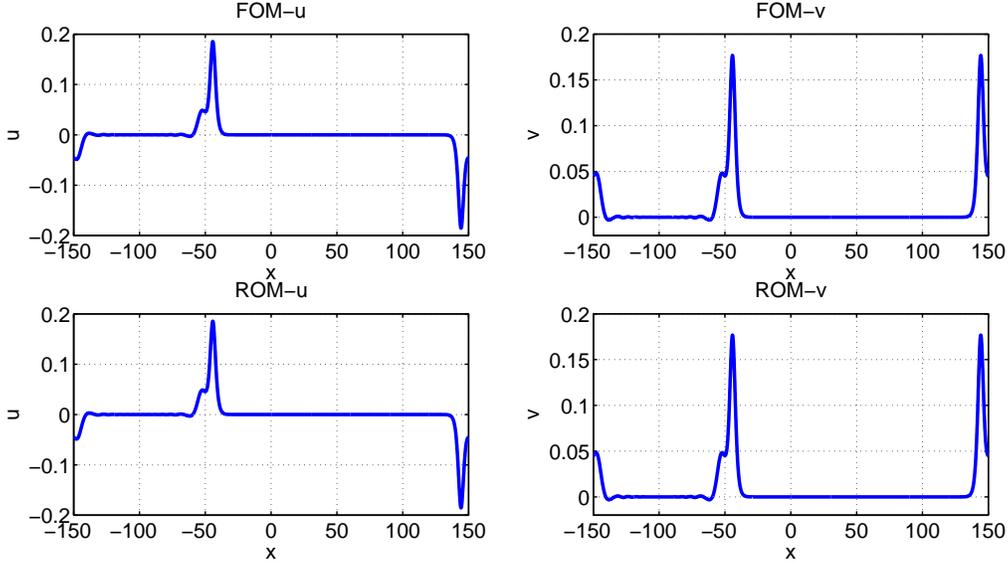}
\caption{FOM and ROM solutions at $T=50$.  }
\label{ex2romsol}
\end{figure}

\begin{figure}[H]
\centering
\includegraphics[width=0.95\linewidth]{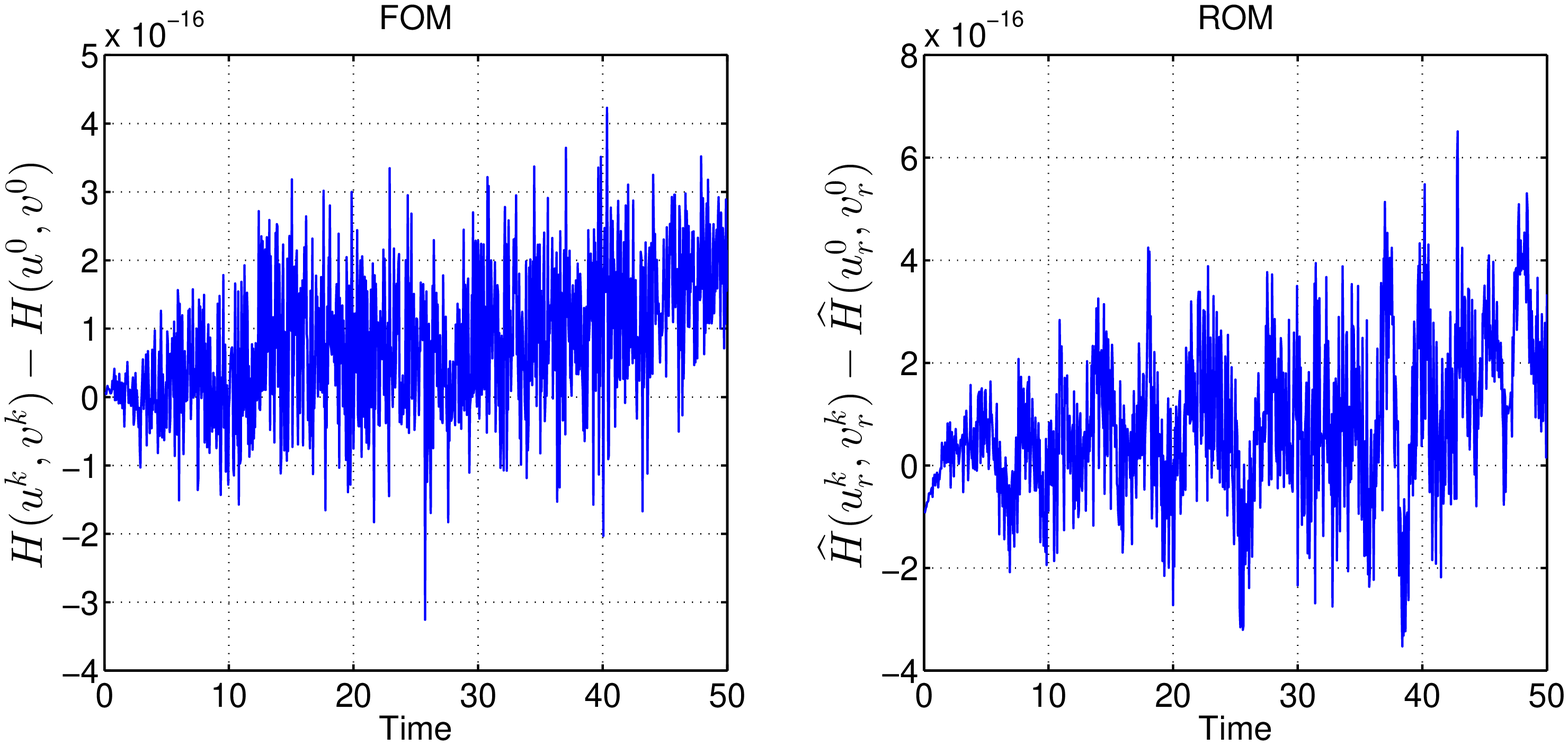}
\caption{Time evolution of the full (left) and the reduced (right) Hamiltonian errors.}
\label{ex2ham}
\end{figure}

In Figure \ref{ex2sing} the singular values decay monotonically without reaching a plateau as for the single KdV equation with $\beta = 10$ in  Figure \ref{ex1sing}.  The number of modes is determined again  by the RIC formula \eqref{ric} as $n=30$ and $n=28$ for $\bm{u}$ and $\bm{v}$ components, respectively. The reduced  and full solutions\footnote{Animations are available as the supplementary material "Ex3\_sol.mp4".} in Figure~\ref{ex2romsol} are visually indistinguishable, and again the discrete Hamiltonian\eqref{dham2} is  preserved accurately by the ROMs in Figure~\ref{ex2ham}.

\subsection{Zakharov-Kuznetsov equation}

We simulate cylindrically symmetric waves  of the Zakharov-Kuznetsov equation \eqref{zk}, that are called as bell-shaped pulses
\cite{Chen11,Iwasaki90,Nishiyama12} with $\alpha =6,\; \mu =1$. The initial condition for two pulses is given by
$$
u(x,y,0) = \sum_{j=1}^2 \frac{c_j}{3} \sum_{m=1}^{10} a_{2m} \left ( \cos \left (2m \text{arccot} \left (\frac{\sqrt{c_j}}{2}r_j\right ) \right ) - 1
\right ),
$$
where $c_1$ and $c_2$ are the velocities of the solitary wave solutions, and $r_i$ is defined by $r_i^2 = (x-x_i)^2+(y-y_i)^2, \;i = 1, 2$. The points
$(x_i, y_i)$ are the location of the peak of $u$. The coefficients  $a_{2m}$ are given in \cite{Nishiyama12}.

Numerical solutions are computed in the rectangular space domain $[0,32]\times[0,32]$ and in the time interval $[0,5]$  using a fine discretization both in space and time,  $\Delta x = \Delta y =0.2286$, $\Delta t = 0.01$, to simulate the waves accurately as in \cite{Chen11,Nishiyama12}. The snapshot matrix is of size $19600\times 500$.

\begin{figure}[H]
\centering
\includegraphics[width=0.5\linewidth]{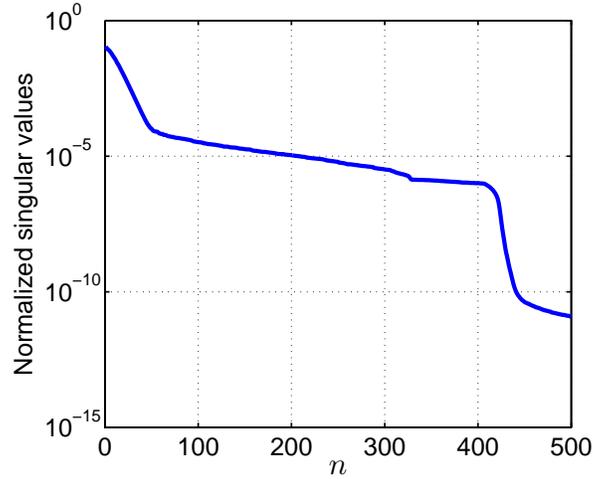}
\caption{Singular values of the snapshot matrix.}
\label{ex3sing}
\end{figure}

The decay of the singular values in Figure~\ref{ex3sing} shows similar behavior as for the single KdV equations in the Figure~\ref{ex1sing} and
in the Figure~\ref{ex2sing}. The number of retained POD modes is $n=50$ according to the RIC formula.

\begin{figure}[H]
\centering
\includegraphics[width=\linewidth] {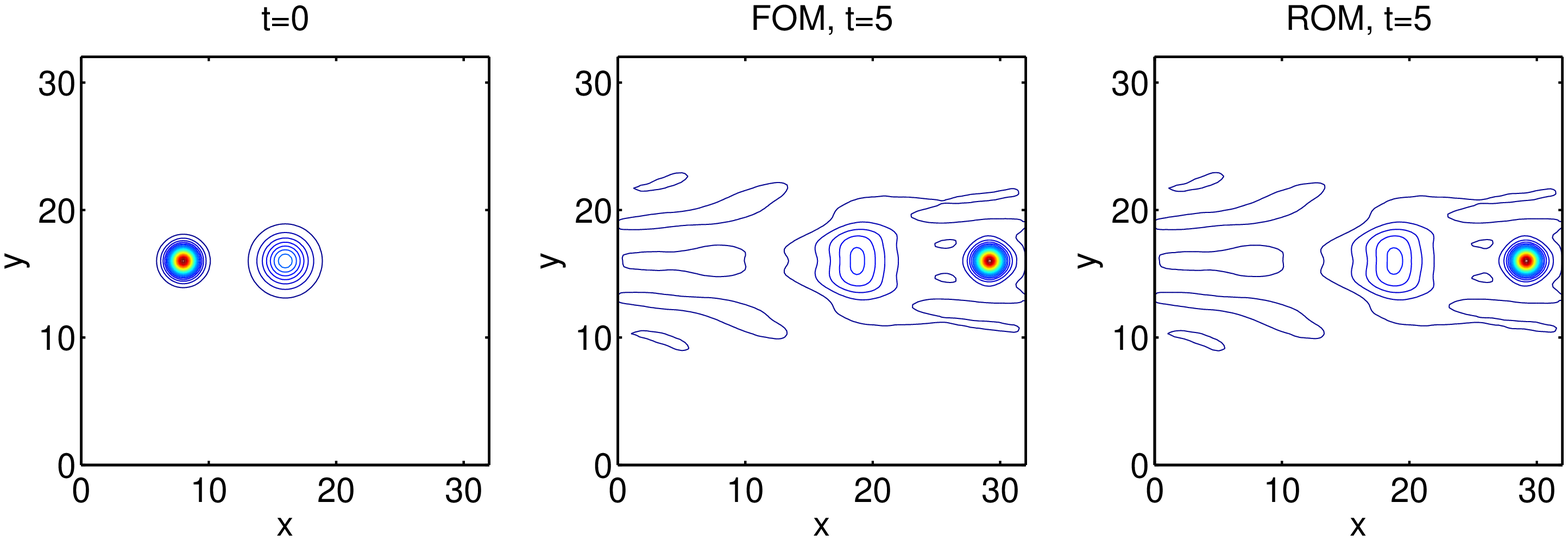}
\includegraphics[width=\linewidth] {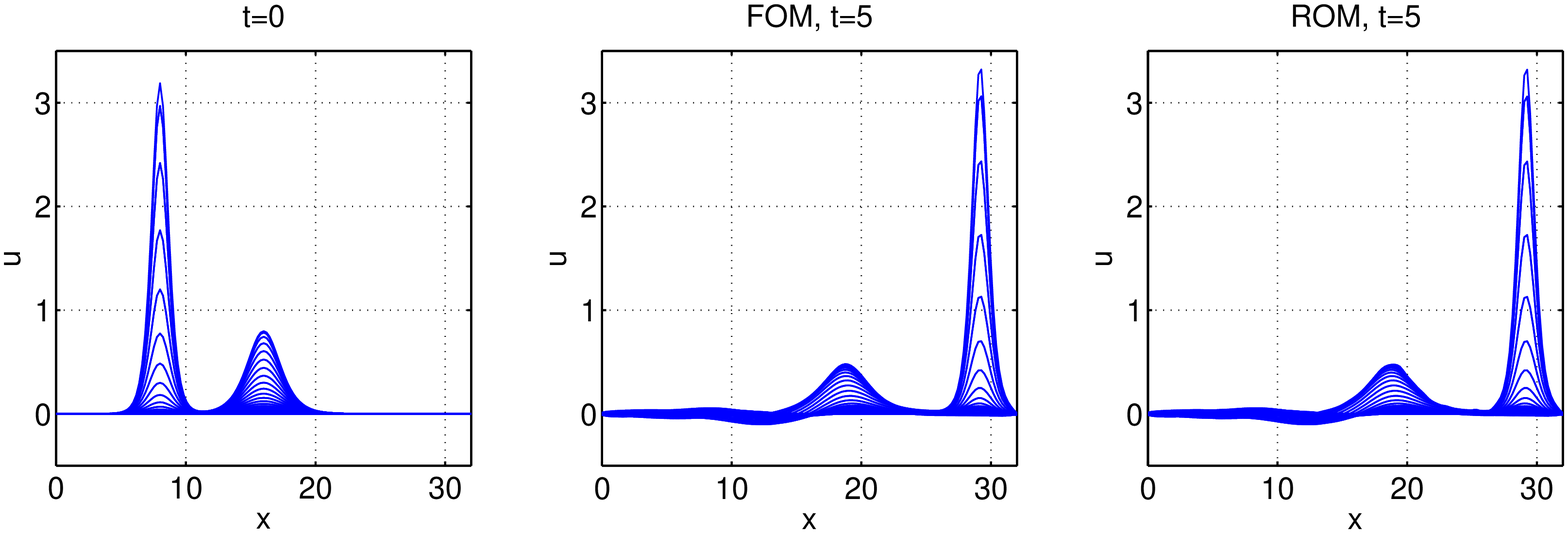}
\caption{Initial profiles,  FOM and  ROM profiles at $t=5$, from left to right.  }
\label{ex3sol}
\end{figure}

\begin{figure}[H]
\centering
\includegraphics[width=0.95\linewidth]{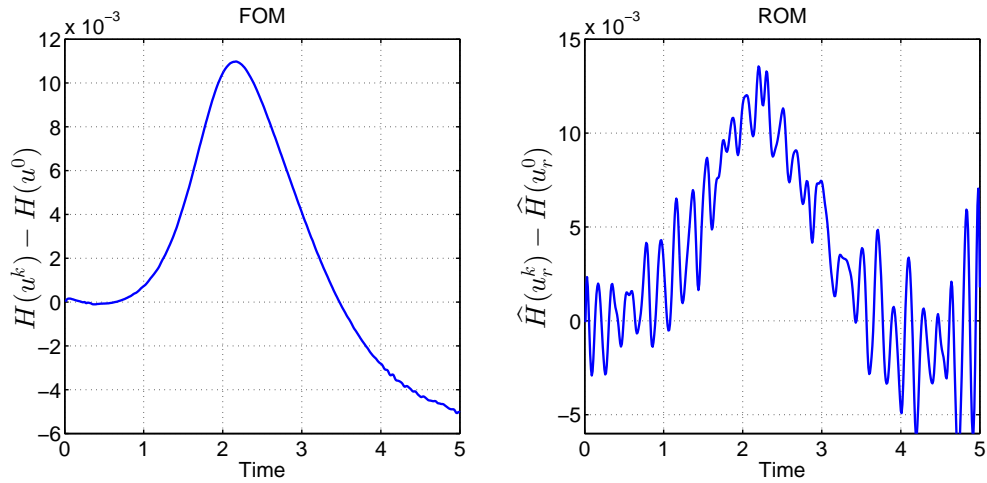}
\caption{Time evolution of the full (left) and the reduced (right) Hamiltonian errors.}
\label{ex3ham}
\end{figure}

In Figure \ref{ex3sol}, the initial profile, the FOM  and ROM profiles\footnote{Animations are available as the supplementary materials "Ex4\_Contour.mp4" and "Ex4\_Piece.mp4".} at the final time $T=5$ are presented. Two dissimilar pulse wave solutions to the Zakharov-Kuznetsov equation at the initial time, evolving in time where the wave structure changes  after
the collision, where the stronger pulse becomes further
stronger and the weaker one gets further weaker after the collision as in  Figure~\ref{ex3sol} by both FOM and ROM.
The discrete Hamiltonian \eqref{ham3d} is well preserved by the ROM in Figure \ref{ex3ham}, even though both Hamiltonian errors are not so small as for the one-dimensional single and coupled KdV equations. However, from the  geometric integration point of view, the Hamiltonian should not drift with time, which is the case for both the full and reduced discrete Hamiltonian in Figure~\ref{ex3ham}.

\subsection{Computational efficiency}

In Table~\ref{cpu}, we present the computational efficiency of the POD and TPOD. The computational cost of the FOM consists of the time required to solve the full solutions, i.e., the creation of snapshots. The computational cost in the offline phase consists of the time required to compute the singular values and singular vectors (POD basis), and the calculation of precomputed matrices. The computational cost in the online phase
consists of the time required to solve the reduced system. In order to measure that to what extend the ROM accelerates the solution process, the speed-up factors are calculated as the ratio of wall-clock time required to solve the FOMs over the wall-clock time required to solve the ROMs in the online phase.
We see that the TPOD approach utilizing MULTIPROD is much faster than the POD, where speed-up factors are given in parenthesis in Table~\ref{cpu}. The efficiency of the TPOD over the POD is much pronounced for the KdV equation with one soliton wave and $\beta =1.5$, and for the Zakharov-Kuznetsov equation, because of larger spatial discretization of the FOMs. In addition, for the single KdV equation with one soliton, the computational efficiency deteriorates with the increasing values of $\beta$ and the number of modes.

\begin{table}[H]
\caption{Wall clock times (in seconds) and speed-up factors (in bold parenthesis) \label{cpu}}
\centering
\resizebox{\textwidth}{!}{
\begin{tabular}{|l|c|c|cc|cc|}
\hline
\multirow{ 2}{*}{System} &  \multirow{ 2}{*}{$n$} & \multirow{ 2}{*}{FOM} &  \multicolumn{2}{c|}{POD}  & \multicolumn{2}{c|}{TPOD} \\
 &  & & Offline & Online &  Offline  & Online \\
	\hline
One soliton ($\beta = 1.5$) & 30 & 178.06 & 5.28 & 62.66 {\bf (2.8)} &  5.67  & 5.55 {\bf (32.1)} \\	
\hline
One soliton ($\beta = 5$) & 60 & 185.62  & 7.49  & 80.70 {\bf (2.3)}  & 8.11    & 46.30 {\bf (4.0)}  \\	
\hline
One soliton ($\beta = 10$) & 90 & 188.81  & 8.25  & 157.34 {\bf (1.2)} &  9.79   &  124.40 {\bf (1.5)} \\	
\hline
Two solitons  & 30 & 7.26  & 1.70  & 2.94 {\bf (2.5)}  &  1.71   & 2.34 {\bf (3.1)}  \\	
\hline
Coupled KdV  & 30, 28 & 17.85  & 1.96  & 2.97 {\bf (6.0)}  & 2.01    & 1.79 {\bf ( 9.9)}  \\	
\hline
Zakharov-Kuznetsov  equation & 50 & 61.15  & 3.33  & 9.25 {\bf (6.6)}  & 3.67  & 0.90 {\bf (68.1)}  \\	
\hline
\end{tabular}}
\end{table}

\section{Conclusions}
\label{con}
We have constructed computationally efficient and accurate ROMs for KdV equations by exploiting the non-canonical Hamiltonian structure.
It is difficult to capture the wave dynamics of PDEs like the KdV equation with a few POD modes. Therefore, in all numerical test problems, the number of the POD modes is relatively large to achieve accurate reduced solutions and to preserve the conserved quantities.
Using TPOD and exploiting the quadratic structure of the KdV equations, the online computational time of ROMs is reduced further.
In a future study, we plan to extend the results of this paper to the parametrized problems using the POD/TPOD-greedy approach in time and in parametric space.\\

\noindent{\bf Acknowledgemets:\/} The authors thank for the constructive comments of the
referees, which helped much to improve the paper.



\end{document}